\documentclass[wide, 10pt]{adonis}

\newtheoremstyle{myplain}
  {}
  {}
  {\itshape}
  {}
  {\sffamily\scshape}
  {.}
  {.5em}
  {}

\theoremstyle{myplain}
\newtheorem{theorem}{Theorem}[section]
\newtheorem{lemma}[theorem]{Lemma}
\newtheorem{proposition}[theorem]{Proposition}
\newtheorem{corollary}[theorem]{Corollary}
\newtheorem*{theorem*}{Theorem}
\newtheorem*{lemma*}{Lemma}
\newtheorem*{proposition*}{Proposition}
\newtheorem*{corollary*}{Corollary}

\newtheoremstyle{mydefinition}
  {}
  {}
  {}
  {}
  {\sffamily\scshape}
  {.}
  {.5em}
  {}

\theoremstyle{mydefinition}
\newtheorem{definition}[theorem]{Definition}
\newtheorem{remark}[theorem]{Remark}
\newtheorem{example}[theorem]{Example}

\newtheorem*{convention*}{Convention}

\newcommand{\diamondend}{%
  \unskip\nobreak\hfill\rotatebox{45}{$\square$}\par
}

\AtEndEnvironment{definition}{\diamondend}
\AtEndEnvironment{remark}{\diamondend}
\AtEndEnvironment{example}{\diamondend}
\AtEndEnvironment{notation}{\diamondend}

\usepackage{bussproofs}
\usepackage{virginialake}
\usepackage{modalops}
\usepackage{quiver}
\usepackage{multicol}
\usepackage{lscape}
\usepackage[c]{esvect}
\usepackage{mathtools}
\usepackage{cancel}

\newcommand\set[1]{\left\{#1\right\}}
\newcommand\domain{\mathrm{Dom}}
\newcommand\aaprox{\approx}
\newcommand\ph{\varphi}

\newcommand\lan[1][]{\mathcal{L}^{\nec}_{#1}}
\newcommand\var{\mathrm{Var}}

\newcommand\necg{\textcolor{gray}{\necf}}
\newcommand\posg{\textcolor{gray}{\posf}}
\newcommand{\Gdl}[1]{\ulcorner #1 \urcorner}
\newcommand{\La}{\mathcal{L}}
\newcommand{\zero}{\overline{0}}

\newcommand\K[1][]{\mathsf{K}_{#1}}
\newcommand\GL{\mathsf{GL}}

\newcommand\ER{\mathsf{ER}}

\newcommand\g[1]{\mathcal{G}#1}

\newcommand\ps[1]{\mathcal{PS}#1}

\newcommand\GLP{\mathsf{GLP}}

\newcommand\prv[1][]{\mathrm{Prv}_{#1}}
\newcommand\prf[1][]{\mathrm{Prf}_{#1}}
\newcommand\con[1][]{\mathrm{Con}_{#1}}

\newcommand\PA{\mathsf{PA}}
\newcommand\EA{\mathsf{EA}}

\newcommand\IS[1]{\mathsf{I}\Sigma_{#1}}

\newcommand\Kax[1][]{\mathrm{K}_{#1}}

\newcommand\NEC[1][]{\mathrm{Nec}_{#1}}

\newcommand\ax{\mathrm{ax}}

\newcommand\botL{{\bot}\mathrm{L}}
\newcommand\botR{{\bot}\mathrm{R}}
\newcommand\toL{{\to}\mathrm{L}}
\newcommand\toR{{\to}\mathrm{R}}

\newcommand\toLsat{{\to}\mathrm{L}^{\mathrm{sat}}}
\newcommand\toRsat{{\to}\mathrm{R}^{\mathrm{sat}}}

\newcommand\modal[2][]{\nec^{#2}_{#1}}

\newcommand\wk{\mathrm{Wk}}
\newcommand\inv[1]{\mathrm{Inv}_{#1}}
\newcommand\ctr{\mathrm{Ctr}}

\newcommand\hg{\mathrm{hg}}
\newcommand\sub{\mathrm{Sub}}

\newcommand\satc[1]{|#1|_{\mathrm{Sat}}}

\newcommand\pI{\texttt{I}} 
\newcommand\pII{\texttt{II}} 

\newcommand\run{\mathscr{R}}

\newcommand\last[1]{\hat{#1}}

\title{Solovay's Construction via Proof Search}
\author{Justus Becker$^{1}$, Borja Sierra Miranda$^{2}$}

\affiliation{$^{1}$University of Birmingham $^{2}$Logic and Theory Group, University of Bern}
\correspondence{$^1$jvb443@student.bham.ac.uk $^2$borja.sierra@unibe.ch}

\runningauthor{J. Becker, B. Sierra Miranda}
\runningtitle{Solovay's Construction via Proof Search}
\abstract{
  Solovay’s celebrated result establishes the arithmetical completeness of \(\GL\) by embedding Kripke models into arithmetic.
  We provide an alternative proof of the arithmetical completeness of \(\GL\), in which we instead embed failed proof searches into arithmetic.
  Our methodology requires only a well-behaved sequent calculus and does not rely on completeness with respect to Kripke semantics or any of its variants.
  This approach is particularly interesting for logics that do not behave well semantically, such as \(\ER\), for which we also provide a new proof of arithmetical completeness.
}

\begin{document}

\maketitle

\section{Introduction}
\emph{Provability logics} are modal logics with an arithmetic interpretation where the modal operators are interpreted as formal provability predicates. 
These logics have a longstanding history, going back to the widely celebrated results by Löb~\cite{lob1955} and Solovay \cite{solovay} of the arithmetic soundness and completeness (respectively) of the modal logic known as Gödel-Löb ($\GL$), characterised by Löb's axiom $\Box (\Box A \to A) \to \Box A$ over the smallest normal modal logic $\K$. 
Since then, this direction of work has been extended to multiple lines of research, including provability under stronger meta-theories, leading to extensions of $\GL$, and the study of multiple provability predicates, be it by virtue of provability in different theories or different notions of provability such as fast provability or $\omega$-provability. 
The former kind of provability logics are also referred to as \emph{modal provability logics} whereas the latter are usually considered in terms of \emph{bimodal provability logics}.

Solovay's proof for \(\GL\) can be divided into two steps.
First, assuming that a formula is not provable in \(\GL\) a \emph{simple} Kripke countermodel of it is obtained.
Second, thanks to the `simplicity' of this model it is shown that it can be embedded inside arithmetic, as desired.
Over time, Solovay's ideas have been proven ubiquitous in the field of provability logic (with a notable exception being~\cite{pakhomov2017solovayscompletenessfixedpoints}).
A few instances of works expanding the ideas of Solovay, among others, are the following.
\begin{itemize}
  \item Carlson~\cite{Carlson1986ModalLW} proved the arithmetical completeness of the bimodal provability logic \(\ER\).
    However, \(\ER\) is Kripke incomplete, so the proof uses a class of models of a Kripke complete sublogic of \(\ER\) which, in addition, globally satisfy enough instances of the \(\ER\)-axiom \(\nec_1(\nec_0 \phi \to \phi)\).
  \item Visser~\cite{Visser1995-VISACO-3} provided an alternative proof of the arithmetical completeness of \(\ER\) using a complete general semantics for \(\ER\), by adding a topology to the Kripke frames.
  \item Japaridze~\cite{japaridze1986modal, japaridze1988polymodal} proved the arithmetical completeness of the multimodal provability logic \(\GLP\).
   As \(\GLP\) is Kripke incomplete, Japaridze used two Kripke models at the same time, each of which preserves some of the structural properties that are necessary for the completeness, but keeping conflicting properties on separate models.
  \item Ignatiev's proof~\cite{Ignatiev_1993} and Beklemishev's proof~\cite{lev-GLP} of the arithmetical completeness of \(\GLP\) both work using models of a Kripke complete sublogic of \(\GLP\) (\(\mathsf{I}\) in the first case and \(\mathsf{J}\) in the second), following the spirit of Carlson's proof for \(\ER\).
\end{itemize}
This state of affairs raises a question: is it possible to avoid the use of semantics, thus avoiding the need of expressing validities of formulas in the logic via structural properties of the relations involved, and still carry out Solovay's construction?
In this paper we will answer positively to this question via the use of sequent calculus and game theory, thus providing a proof that avoids the (explicit) use of Kripke semantics.\footnote{Even if a failed proof search could be viewed as a Kripke model, we will see that this class of models differs from the previous one used in the literature, so it is not immediate from the previous results in the literature that this class can be embedded in arithmetic and provide good witnesses for the unprovability of formulas inside the arithmetic. We will discuss this in greater detail below and in the paper.}

The idea can be sketched as follows.
The essence of Solovay's proof is to embed a \emph{witness} of the unprovability of a formula inside arithmetic.
By a witness of unprovability, we mean an object whose existence guarantees the unprovability of a formula in a fixed modal logic.
Historically, semantics have been used for this purpose since countermodels fulfill this exact role.
On the other hand, proof theory (and, more concretely, the notion of proof) is a witness of the provability of a formula, exactly the opposite of what we want and thus useless for our purpose.
However, game theory comes to the rescue.
A proof can be understood as a winning strategy for one of the two players, usually called verifier, in a game (more specifically, the proof search game).
By the symmetry of games, the absence of a proof can then be understood as the existence of a winning strategy for the other player, refuter.
These winning strategies for refuter will provide us with the necessary witnesses to embed into arithmetic, and in fact it will be of the simplest kind possible: a finite tree.

\textbf{Contributions of the article}.
The contributions of this article are twofold.
\begin{enumerate}
  \item We provide a new methodology for showing arithmetical completeness that, without escaping Solovay's fundamental ideas, avoids the use of modal logic semantics and use instead the existence of a proof search in a well-behaved sequent calculus.
  \item We reprove the arithmetical completeness of \(\GL\) and \(\ER\).
    \(\GL\) will serve as a simple example to understand the methodology.
    However, \(\ER\) will provide a much more interesting example since the proof search approach avoids the complications arising from the semantics and is as simple as the \(\GL\) proof.
    Particularly, the refutations are, in both cases, finite trees and thus simple to embed inside arithmetic.
\end{enumerate}
Additionally, we can see that some uniformity is gained: many provability logics have a simple sequent calculus which is well-behaved with respect to proof search (particularly, this is the case for some Kripke incomplete logics).
We will discuss some possible extensions of the present work in this direction at the conclusion.

\textbf{Relation to semantic proofs}.
As we have explained above, it will not be necessary to use any semantics in our proofs.
However, this does not mean that our approach has no relation to Kripke models.
As is well-known, any failed proof search produces a Kripke model that refutes the sequent and in fact, our work can also be understood from the semantic point of view through an embedding of this model inside arithmetic.
In order to compare the class of models we are embedding into arithmetic with the traditional approaches let us explain a little more about these approaches.

Let us start with the most interesting case, which is the bimodal logic \(\ER\).
The idea is that there is no class of Kripke frames that corresponds to \(\ER\), the reason being that any \(\ER\)-frame would need two irreflexive relations \(R_0, R_1\) such that \(w R_1 v\) implies \(v R_0 v\), among other things.
Then, as \(R_0\) is irreflexive, we would have that \(R_1\) is empty and this would make \(\nec_1 \bot\) a theorem of \(\ER\), but it is not.
We can notice that the tension arises from having two structural properties on the relations that are conflicting in the Kripke setting: the irreflexivity of \(R_0\) and that \(w R_1 v\) implies \(v R_0 v\).
These conditions arise respectively from the axioms \(\nec_0(\nec_0 \phi \to \phi) \to \nec_0 \phi\) and \(\nec_1(\nec_0 \phi \to \phi)\) of \(\ER\).
Each of the proofs of arithmetical completeness of \(\ER\) can be understood as taking one of two natural approaches to solve the tension: either expressing validity in a way that has no relation with the (geometrical) structure of the models or increasing the structural expressive power of the models.
\begin{enumerate}
  \item Carlson's solution is to reduce the number of properties expressed structurally while validating enough instances of axiom \(\nec_1(\nec_0 \phi \to \phi)\) using valuations.
    However, the models used by Carlson still express some validities of formulas structurally such as \(\nec_0 \phi \to \nec_1 \phi\).
  \item Visser's solution is to increase the possibility of expressing validities of formulas structurally by adding a topology in addition to the two relations.
    This increase in expressive power makes the models more complex: they can no longer be assumed to be finite, but just compact with respect to the added topology.
    Still, Visser shows that this suffices for embedding them into arithmetic.
  \item Our solution solves the tension by using proof theoretic methods instead of Kripke semantics, so we do not express \emph{any} validity of axioms structurally.
  In this sense, we follow the same direction as Carlson by expressing less validities structurally, but we take a more radical approach by not expressing any validity structurally.\footnote{Notice that since, in the semantic reading, we produce a finite tree Kripke model the validities of some axioms like axiom (K) or L\"ob's axiom holds in our setting. However, even if this was not the case, the validities of enough instances of these axioms would hold by virtue of the shape of the rules of sequent calculus so we expect this methodology to work in cases where we do not want to have these axioms (see the Future Work section). In fact, producing a finite tree has more to do with producing a \emph{simple} witness to embed inside arithmetic where Solovay's construction is straightfoward than with the validity of axioms of the desired modal logic.}
  Still we can have a semantic reading of our method: we will embed into arithmetic a class of Kripke models that, while witnessing the unprovability in \(\ER\) of the non-theorems, they will not satisfy all the axioms of \(\ER\).
  In fact, a countermodel of \(\phi\) in our setting will only satisfy enough axioms of \(\ER\) to guarantee its unprovability.
\end{enumerate}
We could establish the same analogy for \(\GL\).
So the semantic-minded reader may find interesting our method since we establish the arithmetical completeness of \(\GL\) and \(\ER\) with a different class of models than the known proofs.
We will see, especially in the case of \(\ER\), that this will force some changes from the traditional proof particularly when establishing the Truth Lemma.
More comments about understanding our failed proof search as Kripke models could be found in the subsections devoted to proving the arithmetical completeness of \(\GL\) and \(\ER\), respectively.\footnote{After all these remarks, we want to remind the reader that understanding our method semantically is \emph{possible}, but \emph{not necessary}. It is useful to understand the differences with previous proofs in the literature, but the authors believe that our proofs are much more natural in the proof search framework.}

Finally, let us make a comment about the relation between our work and Japaridze's proof of the arithmetical completeness of \(\GLP\).
The first comment to make is that a failed proof search could be seen as a filtration of the canonical model of the logic, this makes our construction closer to Japaridze's proof than to the above-mentioned constructions.
However, there are still notable differences: Japaridze's proof still uses structural properties of the models he considers.
In fact, to keep all the structural properties that are necessary for Kripke models he separates the properties that are in conflict in two different Kripke models.
Again, the essence of our approach is the use of proof theoretic methods, so the main difference with Japaridze's line of work is the absence of the usual structural properties from the relation that we will embed in arithmetic\footnote{In fact, the differences between our method and Japaridze's method could be even better appreciated in an upcoming work by the second author using proof search to make Solovay's construction in \(\GLP\).}.

\textbf{Proof theory for provability logics}.
Since the (structural) proof theoretic methods will be our main tool, it is worth briefly discussing the literature of proof theory for provability logics.
One of the first instances of structural proof theory for provability logic was Valentini's proof of cut elimination for a sequent calculus for \(\GL\) \cite{Valentini}, which was claimed to contain a gap and multiple authors provided different alternative proofs \cite{borga,mints,sasaki,negri-proof-analysis, negri-proofs-and-countermodels}.
The matter was settled by a formalization of Valentini's proof showing that it is, indeed, correct~\cite{valentini-resolved, valentini-hol}.
Nevertheless, this brief view on the history of the problem shows us that the proofs for cut elimination for \(\GL\) are far from trivial.
This motivated further study of cut elimination in provability logics giving rise to two different approaches: one using a novel measure to reduce the cut size~\cite{Brighton, ian-cut-elim} and another approach keeping the usual measure but allowing non-wellfounded proofs~\cite{coalgebraic} (this last approach inspired by an earlier cut-elimination for non-wellfounded proofs using ultrametric spaces~\cite{shamkanovGrz, shamkanovwGrz}).
These two novel approaches have simplified the treatment of structural proof theory for provability logics as the subsequent works have shown, specifically providing new sequent calculi for provability modal logics which lack one or improving existing ones (see \cite{iSL-cut, iSL-uip, KM-uip, gls-lip, gls-ulip, uip-interpretability, sacchetti-ulip, proofth-interpretability, bimodal-ulip} among others).
With these improvements on the proof theory of provability logic it seems natural to put the sequent calculi to use trying to narrow the gap between semantic and syntactic methods for provability logics.
Such an endeavour is another motivation for this paper, in which we succeed by providing the first Solovay construction replacing semantics with sequent calculus.

\textbf{Structure of the paper.}
In Section~\ref{sec:preliminaries} we will introduce all the concepts necessary to understand the paper, particularly we will talk about arithmetical theories, modal logic, proof theory and game theory.
In Section~\ref{sec:GL} we will describe the methodology of arithmetical completeness via proofs for the first time while applying it to \(\GL\).
In Section~\ref{sec:ER} we will apply our methodology to \(\ER\), while we observe that only minor modifications are needed.
Finally, in Section~\ref{sec:future-work} we will provide some interesting directions for future work with the ideas exposed in the paper.


\section{Preliminaries}\label{sec:preliminaries}

In this section we introduce the necessary notions for the rest of the paper.

\subsection{Modal Logic}
\label{subsec:modal-logic}

We fix an infinite countable set \(\var\) whose elements will be called \emph{propositional variables}.

\begin{definition}
  We define the \emph{language of modal logic} \(\lan[1]\) by the following grammar
  \[
    \phi ::= p \mid \bot \mid \phi \to \phi \mid \nec_0 \phi,
  \]
  where \(p \in \var\) and the \emph{language of bimodal logic} \(\lan[2]\) by the following grammar
  \[
    \phi ::= p \mid \bot \mid \phi \to \phi \mid \nec_0 \phi \mid \nec_1 \phi,
  \]
  where \(p \in \var\).
  The expressions of \(\lan[1]\) are called \emph{modal formulas} and the expressions of \(\lan[2]\) are called \emph{bimodal formulas}.
\end{definition}
We will write \(\phi = \psi\) to mean that \(\phi\) and \(\psi\) are syntactically equal.
As usual, we define the following abbreviations:
\begin{align*}
  &\neg \phi := \phi \to \bot,
  &&\phi \vee \psi := \neg \phi \to \psi,
  &\phi \wedge \psi := \neg(\neg \phi \vee \neg \psi),
  &&\pos_i \phi := \neg \nec_i \neg \phi.
\end{align*}
In case we work with formulas in \(\lan[1]\) we will write \(\nec_0\) simply as \(\nec\).
Given a formula \(\phi\) we define its complexity, denoted \(|\phi|\) as
\begin{align*}
  &|p| = |\bot| = 1,
  &&|\phi_0 \to \phi_1| = |\phi_0| + |\phi_1| + 1,
  &|\nec_i \phi_0| = |\phi_0|+ 1.
\end{align*}
Formulas of shape \(\nec_i \phi\) for \(i \leq 1\) will be called \(\nec\)-formulas, the ones of shape \(\nec_0 \phi\) will be called \(\nec_0\)-formulas and the ones of shape \(\nec_1 \phi\) will be called \(\nec_1\)-formulas.
Notice that in \(\lan[1]\), \(\nec_0\)-formulas and \(\nec\)-formulas are the same (since there are no \(\nec_1\)-formulas in that language).

In this paper we will be working with one modal logic, \(\GL\), and one bimodal logic, \(\ER\).
We now introduce these logics.

\(\GL\) is defined as the smallest set of formulas in \(\lan[1]\) such that
\begin{enumerate}
  \item all propositional tautologies in the language \(\lan[1]\) belong to \(\GL\);
  \item \(\nec(\phi \to \psi) \to \nec \phi \to \nec \psi\) belongs to \(\GL\) for \(\phi, \psi \in \lan[1]\);
  \item \(\nec(\nec \phi \to \phi) \to \nec \phi\) belongs to \(\GL\) for \(\phi \in \lan[1]\);
  \item \(\phi, \phi \to \psi \in \GL\) implies \(\psi \in \GL\);
  \item \(\phi \in \GL\) implies \(\nec \phi \in \GL\).
\end{enumerate}
The formulas in \(\GL\) are called the \emph{theorems of \(\GL\)}.
We will write \(\GL \vdash \phi\) to mean that \(\phi\) is a theorem of \(\GL\).
It is well-known (e.g. \cite[p.~11]{boolos} Theorem 18) that we have the following lemma.
\begin{lemma}
  For any formula \(\phi\), \(\GL \vdash \nec\phi \to \nec \nec \phi\).
\end{lemma}

\(\ER\) is defined as the smallest set of formulas in \(\lan[2]\) such that
\begin{enumerate}
  \item all propositional tautologies in the language \(\lan[2]\) belong to \(\ER\);
  \item \(\nec_i(\phi \to \psi) \to \nec_i \phi \to \nec_i \psi\) belongs to \(\ER\) for \(\phi, \psi \in \lan[2]\), \(i \leq 1\);
  \item \(\nec_i(\nec_i \phi \to \phi) \to \nec_i \phi\) belongs to \(\ER\) for \(\phi \in \lan[2]\), \(i \leq 1\);
  \item \(\nec_i\phi \to \nec_{\overline{i}}\nec_i \phi\) belongs to \(\ER\) for \(\phi \in \lan[2]\), \(i \leq 1\);
  \item \(\nec_1(\nec_0 \phi \to \phi)\) belongs to \(\ER\) for \(\phi \in \lan[2]\);
  \item \(\phi, \phi \to \psi \in \ER\) implies \(\psi \in \ER\);
  \item \(\phi \in \ER\) implies \(\nec_i \phi \in \ER\) for \(i \leq 1\).
\end{enumerate}
The formulas in \(\ER\) are called the \emph{theorems of \(\ER\)}.
We will write \(\ER \vdash \phi\) to mean that \(\phi\) is a theorem of \(\ER\).
Using the same reasoning as for \(\GL\) we have the following lemma.
\begin{lemma}
  For any formula \(\phi\) and \(i \leq 1\), \(\ER \vdash \nec_i\phi \to \nec_i \nec_i \phi\).
\end{lemma}

Let us conclude this section, even if it is not strictly related, by recalling the definition of finite rooted tree.
A \emph{finite rooted tree} is a pair \((W,R)\) such that
\begin{enumerate}
  \item \(W\) is a finite set and \(R\) is a relation on \(W\),
  \item \(R\) is a transitive, irreflexive relation
  \item There is a unique node \(r\) called \emph{the root} such that for any \(w \in W \setminus \set{r}\), \(r \mathrel{R} w\).
\end{enumerate}
Usually, \(w R v\) is read as \(v\) is a \emph{successor} of \(w\) or \(w\) is an \emph{ancestor} of \(v\).
Additionally, given a finite tree \((W,R)\) we define the \emph{immediate successor} relation as \(w R^{1} v\) if \(w R v\) and there is no \(u \in W\) such that \(w R u R v\).
If \((R^{1})^+\) denotes the transitive closure of \(R^{1}\) we note that \((R^{1})^+ = R\).
This \emph{immediate successor} relation will be more important for our purposes than the relation \(R\).
Finally, we remember a basic property of finite trees: if \(w,v \in W\) such that \(w \neq v\), \(w \mathrel{\not R} v\) and \(v \mathrel{\not R} w\) then there exist nodes \(u,u_0,u_1\) such that \(u_0,u_1\) are immediate successors of \(u\), (\(u_0 = w\) or \(u_0 R w\)) and (\(u_1 = v\) or \(u_1 R v\)).
In words, if \(w\) and \(v\) do not belong to the same branch they must share an ancestor where the branch of \(w\) and the branch of \(v\) are divided.

\subsection{Arithmetical theories}

We assume the reader is familiar with first order logic and theories of arithmetic.
Solovay’s Theorem, and its variants, apply to a very large class of theories. 
We wish to make this fact as explicit as possible by stating Solovay-like theorems in their most general form, removing unnecessary details and allowing us to see what properties we need to require from theories, such that Solovay’s proof strategy applies to them.
We thus avoid talking about specific theories as much as possible, and instead only list the essential properties that will be necessary for our proof.

For simplicity we will consider theories in the \emph{signature of arithmetic with exponentiation} which has functions $(\zero , s, + , \times, \mathrm{exp})$ of arities \((0,1,2,2,2)\) and a binary relation symbol \(\prec\).
We use  $\approx$, \(\prec\), and $\zero$ to distinguish equality, the natural numbers order, and numbers in the object language from the same notions in the meta language.

We define the symbol $\precsim$ as the formula with two free variables:
\[
x \precsim y := x \approx y \vee x \prec y. 
\]
We call a formula $\ph$ of $\La$ \emph{bounded} if every quantifier is guarded as $\exists x \prec t \, \psi$  or $\forall x \prec t \, \psi$, respectively short for $\exists x ( x \prec t \wedge \psi)$ and $\forall x ( x \prec t \to \psi)$, where $t$ is some term of in $\La$ in which \(x\) does not occur.
We call the set of bounded formulas $\Delta_0 = \Sigma_0 = \Pi_0$ and define the arithmetic hierarchy as usual: $\Sigma_{n+1} := \{ \exists x_1 ... \exists x_n \ph \mid \ph \in \Pi_n\}$ and $\Pi_{n+1} := \{ \forall x_1 ... \forall x_n \ph \mid \ph \in \Sigma_n\}$ with $n \in \omega$.
We will say that a formula is \emph{strictly \(\Sigma_n\)} or \emph{strictly \(\Pi_n\)} if each block of \(\forall\) and \(\exists\) quantifiers in the previous definition consists of exactly one quantifier.

We inductively define \emph{numerals} $\overline n$ as terms in $\La$ with $\overline{n+1} = s(\overline{n})$ and $\zero$ being the zero constant.
Let $\overline{\omega}$ denote the set of numerals.
We call a \emph{Gödel coding} for the language of $T$ an injective function $\Gdl{\cdot} : \La \to \overline{\omega}$ from formulas of $T$ to numerals.
Additionally, we will say that
\begin{enumerate}
  \item  \emph{satisfies uniform diagonalization in \(T\)} if for any $\vec x = (x_1 , ... ,x_n)$ tuple of variables and $P_1 (\vec x) , ... , P_n (\vec x )$ formulas in the language of \(T\), there are sentences $S_1 , ... , S_n$ such that 
    \[
      T \vdash S_i \leftrightarrow P_i (\Gdl{S_1}, \ldots , \Gdl{S_n}).
    \]
  for all $i \leq n$,
  \item \emph{has \(\Delta_0\)-negation} if there a \(\Delta_0\)-formula \(\mathrm{Neg}(x,y)\) such that \(T \vdash \mathrm{Neg}(\overline{n},\overline{m})\) if and only if \(\overline{n} = \Gdl{\phi}\) and \(\overline{m} = \Gdl{\neg \phi}\) for some formula \(\phi\).
\end{enumerate}

Among the theories we are going to consider the theory of \emph{elementary arithmetic}, denoted \(\EA\), will be of special importance.
This is the theory with the defining axioms of exponentiation and induction for bounded formulas (see~\cite{HajekPudlak1998} for details).
We remember some important facts about \(\EA\).

\begin{proposition}\label{properties-EA}
  We have the following.
  \begin{enumerate}
    \item \(\EA\) proves \(\Sigma_0\)-minimization, i.e., for any \(\Sigma_0\)-formula \(\phi(x)\) (possibly with more free variables) we have that the universal closure of \( \exists x (\phi(x)) \to \exists x(\phi(x) \wedge \forall y \prec x (\neg \phi(y))) \) is provable in \(\EA\).
    \item  \(\EA\) is \(\Sigma_1\)-complete, i.e., for any \(\Sigma_1\)-sentence \(\alpha\) we have that \(\mathbb{N} \vDash \alpha\) implies \(T \vdash \alpha\).
    \item  In \(\EA\) every \(\Sigma_1\)-formula is equivalent to a strict \(\Sigma_1\)-formula.
  \end{enumerate}
\end{proposition}

With \(\EA\) at hand we are ready to define what will be an arithmetic theory in this paper.

\begin{definition}
  An arithmetic theory \(T\) will be a consistent first order theory whose signature is the signature of arithmetic with exponentiation and the following properties.
  \begin{enumerate}
    \item it extends \(\EA\), i.e.,\ if \(\EA \vdash \alpha\) then \(T \vdash \alpha\).
    \item there is a G\"odel coding \(\Gdl{\cdot}\) satisfying uniform diagonalization in \(T\) and with \(\Delta_0\)-negation.
  \end{enumerate}
Some well-known examples of arithmetical theories are \(\PA\), i.e., Peano Arithmetic, and multiple of its fragments such as the theories \(\IS{n}\) for \(n \geq 1\) (see~\cite{HajekPudlak1998} for details).
\end{definition}
Note that by virtue of extending \(\EA\), arithmetical theories will also satisfy the properties of \(\EA\) at Proposition~\ref{properties-EA}.





It will also be handy, given strict \(\Sigma_1\)-formulas \(\phi = \exists x(\phi_0(x))\), \(\psi = \exists y(\psi_0(y))\) to define the following formulas:
\[
  \phi < \psi := \exists x (\phi_0(x) \wedge \forall y \precsim x (\neg \psi_0(y))),
  \qquad
  \phi \leq \psi := \exists x (\phi_0(x) \wedge \forall y \prec x (\neg \psi_0(y))).
\]
We note that, by definition, \(\phi < \psi\) and \(\phi \leq \psi\) are \(\Sigma_1\)-formulas.

Finally, let us see how to interpret modal formulas as arithmetical sentences.
\begin{definition}
	Fix some arithmetical predicates $P(x), Q(x)$ with exactly one free variable.
	We call a map $(\cdot)^*_P: \lan[1] \to \La$ an \emph{arithmetic interpretation} if $(\nec \phi)^\ast = P(\Gdl{\phi^\ast})$, $\bot^\ast = (\zero \aaprox \overline 1)$ and $(\cdot)^*$ commutes with Boolean connectives. 

	We call a map $(\cdot)^*_{P,Q}: \lan[2] \to \La$ an \emph{arithmetic interpretation} of the pair $(P, Q)$ if $(\nec_0 \phi)^* = P(\Gdl{\phi^*})$, $(\nec_1 \phi)^\ast = Q(\Gdl{\phi^*})$, $\bot^* = (\zero \aaprox \overline 1)$ and $(\cdot)^*$ commutes with Boolean connectives. 
\end{definition}

\subsection{Proof Theory}

Let us start by fixing the notions concerning proof theoretic trees.\footnote{We call \emph{proof theoretic} trees the structure that we are going to use to present proofs, in contrast with the notion of tree at the end of Section~\ref{subsec:modal-logic}.}
Given a set \(X\) we will write \(X^*\) (\(X^+\)) to mean the set of (non-empty) finite sequences of elements in \(X\), \(\epsilon\) will denote the empty sequence.
As usual, we will write \(w \leq v\) to mean that \(w\) is an initial prefix of \(v\) and \((w,v]\) to mean the set \(\set{u \in X^* \mid w < u \leq v}\).
A (finitely branching) proof theoretic tree on \(A\) is a function \(T\) whose image is contained in \(A\) and whose domain is a prefix-closed non-empty subset of \(\mathbb{N}^*\) such that for every \(w \in \domain(T)\) there is a unique natural number \(k\) (called the \emph{arity of \(w\) in \(T\)}) such that \(wi \in \domain(T)\) iff \(i < k\).
The \(wi \in \domain(T)\) are also called the \emph{immediate successors of \(w\)}.
Note that the domain of any proof theoretic tree always contains the word \(\epsilon\), and is called the \emph{root of \(T\)}.

The elements of \(\domain(T)\) are also called the \emph{nodes of \(T\)}, the \(0\)-ary nodes are called \emph{leaves} and the rest of the nodes are called \emph{interior nodes}.
An \emph{infinite branch in a proof theoretic tree \(T\)} is sequence of nodes \((w_i)_{i \in \mathbb{N}}\) such that \(w_0 = \epsilon\) and for each \(i\) there is a \(j\) such that \(w_{i+1} = w_i j\).
Finally, a proof theoretic tree is said to be \emph{wellfounded} if it contains no infinite branches.

We fix a set \(\text{Seq}\) whose elements we will call \emph{sequents}, a sequent rule is a subset of \(\text{Seq}^+\).
Given a rule \(R\), we say that \((S_0,\ldots,S_{n-1},S)\) is an \emph{instance of \(R\) with premises \(S_0,\ldots,S_{n-1}\) and conclusion \(S\)} if \((S_0,\ldots,S_{n-1},S) \in R\).
A rule \(R\) is said to be \emph{\(n\)-ary} if \(R \subseteq \text{Seq}^{n+1}\).
It will be common to write an instance \((S_0,\ldots,S_{n-1},S)\) of a rule \(R\) as
\[
  \AxiomC{\(S_0\)}
  \AxiomC{\(\cdots\)}
  \AxiomC{\(S_{n-1}\)}
  \RightLabel{\(R\)}
  \TrinaryInfC{\(S\)}
  \DisplayProof
\]
We introduce what we mean precisely with a sequent calculus.

\begin{definition}
  A \emph{sequent calculus} is a set of sequent rules. 
  Given a sequent calculus \(\mathcal{G}\), a \emph{proof in \(\mathcal{G}\)} is a wellfounded proof theoretic tree with labels in \(\text{Seq} \times \mathcal{G}\) such that for each node \(w\) with immediate successors \(w0,\ldots,w(n-1)\) we have that \((S_0,\ldots,S_{n-1},S) \in R\) where \(S\) is the sequent at \(w\), \(R \in \mathcal{G}\) is the rule at \(w\) and \(S_i\) is the sequent at \(wi\).
\end{definition}

As usual, we will write \(\pi \vdash_{ \mathcal{G}} S\) to mean that \(\pi\) is a proof in \(\mathcal{G}\) whose root is annotated with the sequent \(S\) and \(\mathcal{G}\vdash S\) to mean that there is a \(\pi\) such that \(\pi \vdash_{\mathcal{G}} S\).
We may drop the subscript \(_{\mathcal{G}}\) from the first notation in case it can be inferred from context.

Rules can interact with sequent calculi in different ways, the following definition introduces some of these interactions.

\begin{definition}\label{def:prop-of-rules}
  Let \(\mathcal{G}\) be a sequent calculus and \(R\) be a rule.
  We say that \(R\) is
  \begin{enumerate}
    \item \(R\) admissible if for any instance \((S_0,\ldots,S_{n-1},S) \in R\) we have that \(\pi_i \vdash_{\mathcal{G}} S_i\) for \(i < n\) implies the existence of \(\pi \vdash_{\mathcal{G}} S\).
      In addition we say that \(R\) is \emph{admissible preserving height} if \(\hg(\pi) \leq \max_{i < n}\hg(\pi_i)\).
    \item An \(n\)-ary rule \(R\) is \(i\)-invertible  in \(\mathcal{G}\) (where \(i < n\)) if the rule
      \(
        \set{(S,S_i) \mid (S_0,\ldots,S_{n-1},S) \in R}
      \)
      is admissible.
      \(R\) is \emph{invertible} if it is \(i\)-invertible for each \(i < n\).
  \end{enumerate}
We will talk about invertibility preserving height with the obvious meaning.
\end{definition}
Given an admissible rule \(R\) in \(\mathcal{G}\), an instance \((S_0,\ldots,S_{n-1},S) \in R\) and proofs \(\pi_i \vdash_{\mathcal{G}} S_i\) for \(i < n\), we will write \(R(\pi_0,\ldots,\pi_{n-1})\) to mean a proof of \(S\) in \(\mathcal{G}\) that exists by admissibility of \(R\).
For the invertibility of a rule \(R\) we will write instead \(\inv{R}(\pi)\).

\subsection{Game Theory}
We need some basic notions and results in combinatorial game theory.
We start with the definition of game.

\begin{definition}
  A \emph{ruleset} is a tuple \(\mathscr{G} = (P,\leadsto, P_{\pI})\) such that
  \begin{enumerate}
    \item \(P\) is a set whose elements are called \emph{positions},
    \item \(\leadsto\) is a binary relation on \(P\), for \(G \in P\) the elements of \(\set{G' \mid G \leadsto G'}\) are called the \emph{moves at \(G\)},
    \item \(P_{\texttt{I}} \subseteq P\) whose elements are called \emph{positions of player \(\pI\)}, the elements of \(P_{\pII} = P \setminus P_{\pI}\) are called \emph{positions of player \(\pII\)}.
  \end{enumerate}
  A \emph{game} is a pair \((\mathscr{G},G)\) where \(\mathscr{G}\) is a ruleset and \(G\) is a position of \(\mathscr{G}\).
  Such a pair will usually be denoted as \(\mathscr{G} @ G\) or just \(G\) if \(\mathscr{G}\) is clear from context.
\end{definition}

In the following definition we explain how games are played and won.

\begin{definition}
  Given a game \(\mathscr{G} @ G\) we define:
  \begin{enumerate}
    \item A \emph{finite run} of \(\mathscr{G} @ G\) is a non-empty sequence \((G_i)_{i < k}\) where \(k \in \mathbb{N} \setminus \set{0} \) such that \(G_0 = G\) and \(G_i \leadsto G_{i+1}\) for \(i < k\).
      The set of finite runs of \(\mathscr{G} @ G\) is denoted as \(\run_{<\omega}(\mathscr{G} @ G)\) or \(\run_{<\omega}(G)\) if the ruleset is clear from context.
    \item An \emph{infinite run} of \(\mathscr{G} @ G\) is a sequence \((G_i)_{i < \omega}\) such that \(G_0 = G\) and \(G_i \leadsto G_{i+1}\) for \(i < \omega\).
      The set of infinite runs of \(\mathscr{G} @ G\) is denoted as \(\run_{\omega}(\mathscr{G} @ G)\) or \(\run_{\omega}(G)\) if the ruleset is clear from context.
    \item A \emph{run} of \(\mathscr{G} @ G\) is a finite or infinite run of \(\mathscr{G} @ G\).
      The set of infinite runs of \(\mathscr{G} @ G\) is denoted as \(\run(\mathscr{G} @ G)\) or \(\run(G)\) if the ruleset is clear from context.
    \item A \emph{play} of \(\mathscr{G} @ G\) is just a maximal run of \(\mathscr{G} @ G\), i.e., it is a run \((G_i)_{i < K}\) such that either \(K = \omega\) or \(K \in \mathbb{N} \setminus \set{0}\) and there are no possible moves at \(G_{K-1}\).
    A finite play is \emph{won by \(\pI\)} (resp.\ \(\pII\)) if its last element belongs to \(P_\pII\) (resp.\ \(P_\pI\)).\qedhere
  \end{enumerate}
  Given a run \(w = (G_0,\ldots,G_{k-1})\) we will write \(\hat{w} = G_{k-1}\), i.e., \(\hat{w}\) denotes the last position of the run.
\end{definition}

Finally, we can introduce the notion of strategy.
A strategy is a function that tells a player which move to play in each position, based on the history of the game.

\begin{definition}
  A strategy for \(\pI\) in a game \(\mathscr{G}@G\) is a function \(s : \set{w \in \run_{< \omega}(G) \mid \hat{w} \in P_{\pI}} \longrightarrow P\) such that for any finite run \(w\) we have that \(\hat{w} \leadsto s(w)\).
  Strategies for \(\pII\) are defined similarly.

  For a strategy \(s_{\pI}\) for \(\pI\) and a run \((G_i)_{i < k}\) of \(\mathscr{G} @ G\), we say that \((G_i)_{i < k}\) is \emph{\(s_{\pI}\)-consistent} if for any \(i < k - 1\) such that \(G_i \in P_{\pI}\) we have that \(G_{i+1} = s(G_0, \ldots, G_i)\).
  Consistent runs for strategies of \(\pII\) are defined similarly.
  For a strategy \(s\) (of player \(\pI\) or \(\pII\)), the set of \(s\)-consistent (finite, infinite) runs is denoted as \(\run(G | s)\) (\(\run_{<\omega}(G | s)\), \(\run_{\omega}(G | s)\)).

  A strategy \(s\) for \(\pI\) is said to be a \emph{winning strategy in the game \(\mathscr{G} @ G\)} if every \(s\)-consistent play of \(\mathscr{G} @ G\) is won by \(\pI\).
  Winning strategies for \(\pII\) are defined similarly.
\end{definition}

Finally, we want to introduce a class of well-behaved games.
All the games we are going to consider belong to this class.

\begin{definition}
  Given two games \(G, G'\) we say that \(G'\) is a \emph{subgame} of \(G\) if \(G = G'\) or \(G \leadsto^+ G'\), where \(\leadsto^+\) is the transitive closure of \(\leadsto\).

A game \(\mathscr{G} @ G\) is called \emph{short} if
\begin{itemize}
  \item \(\mathscr{G} @ G\) has finitely many different subgames, and
  \item \(\mathscr{G} @ G\) has no infinite runs.
\end{itemize}
\end{definition}

Note that the second condition implies that \(\leadsto\) on the subgames of \(G\) is conversely wellfounded and together with the first condition it is a finite conversely wellfounded relation on the subgames of \(G\).
This implies that for a short game \(\mathscr{G}@G\) the set \((\run(G), \leadsto^+)\) of runs with the strict subgame relation is a finite tree whose immediate successor relation is \(\leadsto\).
The trees we will embed in arithmetic will be found inside this tree.

We conclude the section with a fundamental theorem of combinatorial game theory.
\begin{theorem}[Determinacy of short games~\cite{Zermelo1913}]
  For any short game exactly one of the players has a winning strategy.
\end{theorem}

\section{Arithmetical completeness for \(\GL\)}\label{sec:GL}

We are almost ready to begin the proof of arithmetical completeness for \(\GL\).
First, we still need one crucial notion: that of a Hilbert-Bernays provability predicate.

\begin{definition}
  Let \(T,B\) be arithmetical theories.
  A \emph{(single conclusion) Hilbert-Bernays provability predicate \(\prv(x)\) for \(T\) in \(B\)} is an arithmetical formula with exactly one free variable such that
  \begin{enumerate}
    \item \(\prv(x)\) is a strict \(\Sigma_1\)-formula of shape \(\exists p.\ \prf(p,x)\), for any arithmetical sentence \(\alpha\) let us denote \(\prv(\Gdl{\alpha})\) as \(\necg \alpha\);
    \item (\(\NEC[\necg]\)) \(T \vdash \alpha\) implies \(B \vdash \necg \alpha\) for any arithmetical sentence \(\alpha\),
    \item (\(\Kax[\necg]\)) \(B \vdash \necg(\alpha \to \beta) \to \necg \alpha \to \necg \beta\) for any arithmetical sentences \(\alpha, \beta\),
    \item (\(\Sigma_1\)-completeness) \(B \vdash \alpha \to \necg \alpha\) for any \(\Sigma_1\)-sentence \(\alpha\),
    \item (Single conclusion) \(B \vdash \forall x,y.\ \prf(p,x) \wedge \prf(p,y) \to x \approx y\), and
  \end{enumerate}
  When \(T = B\) we will simply say that \(\prv(x)\) is a \emph{Hilbert-Bernays provability predicate for \(T\)}.
  \(\prf(p,x)\) is called a \emph{proof predicate}, more concretely because of the last condition we say that it is a \emph{single-conclusion proof predicate}.
  Additionally, we will say that \(\prv(x)\) is \emph{sound} if for any sentence \(\alpha\), \(\mathbb{N} \vDash \alpha\) implies \(T \vdash \alpha\).
\end{definition}

In the rest of the section we are going to show the arithmetical completeness of \(\GL\), more precisely we will be showing the following theorem.

\begin{theorem}[Arithmetical completeness for \(\GL\)]\label{arithmetical-completeness-GL}
  Let \(T\) be a sound arithmetical theory and \(\prv\) be a sound Hilbert-Bernays provability predicate for \(T\).
  Then, for any modal formula \(\phi\) such that \(\GL \nvdash \phi\) there is an arithmetical interpretation \((\cdot)^*_{\prv}\) with \(T \nvdash \phi^*\).
\end{theorem}

For simplicity we have assumed that \(T\) and \(\prv(x)\) are sound.
We have also taken \(T = B\), this are non-essential restrictions and we will discuss some known generalizations at the end of the section.

\subsection{Proofs and refutations}

\begin{figure}
  \[
    \AxiomC{}
    \RightLabel{\(\ax\)}
    \UnaryInfC{\(p, \Gamma \Rightarrow p, \Delta\)}
    \DisplayProof
    \qquad
    \AxiomC{}
    \RightLabel{\(\botL\)}
    \UnaryInfC{\(\bot, \Gamma \Rightarrow \Delta\)}
    \DisplayProof
    \qquad
    \AxiomC{\( \Gamma \Rightarrow \Delta\)}
    \RightLabel{\(\botR\)}
    \UnaryInfC{\( \Gamma \Rightarrow \bot,\Delta\)}
    \DisplayProof
  \]

  \[
    \AxiomC{\(\Gamma \Rightarrow \phi, \Delta\)}
    \AxiomC{\(\psi, \Gamma \Rightarrow \Delta\)}
    \RightLabel{\(\toL\)}
    \BinaryInfC{\(\phi \to \psi, \Gamma \Rightarrow \Delta\)}
    \DisplayProof
    \qquad
    \AxiomC{\(\phi, \Gamma \Rightarrow \psi, \Delta\)}
    \RightLabel{\(\toL\)}
    \UnaryInfC{\(\Gamma \Rightarrow \phi \to \psi, \Delta\)}
    \DisplayProof
  \]

  \[
    \AxiomC{\(\necd \Sigma, \nec \phi\Rightarrow \phi\)}
    \RightLabel{\(\modal{\GL}\)}
    \UnaryInfC{\(\nec \Sigma, \Gamma \Rightarrow \nec \phi, \Delta\)}
    \DisplayProof
  \]
  \caption{\(\GL\) sequent rules}
  \label{fig:GL-sequent-rules}
  \separator
\end{figure}

We start by recalling the well-known sequent calculus for \(\GL\)~\cite{Valentini}.
The rules for this sequent calculus can be found at Figure~\ref{fig:GL-sequent-rules}.
As usual, the displayed formula at the conclusion of \((\botR)\), \((\toL)\), \((\toR)\), \((\modal{\GL})\) is called the \emph{principal formula} of the rule instance.
The formulas in \(\nec \Sigma\) at \((\modal{\GL})\) are called the \emph{auxiliary formulas} of the rule instance.
Together, principal and auxiliary formulas are called the \emph{active formulas} of the instance.
Additionally, in \((\ax)\), \((\botL)\) and \((\modal{\GL})\) the formulas at \(\Gamma\) and \(\Delta\) are called \emph{weakening formulas}.
The sequents which are instances of \((\ax)\) and \((\botL)\) are called \emph{axiomatic sequents}.
An \emph{extended axiomatic sequent} is a sequent \(\Gamma \Rightarrow \Delta\) such that \(\Gamma \cap \Delta \neq \varnothing\) or \(\bot \in \Gamma\).

\begin{definition}
  We define the sequent calculus \(\g{\GL}\) as the one given by the rules of Figure~\ref{fig:GL-sequent-rules}.
\end{definition}

It is straightforward to show that extended axiomatic sequents are provable.
\begin{proposition}
  Let \(S\) be an extended axiomatic sequent, then \(\g{\GL} \vdash S\).
\end{proposition}

It is well-known that \(\g{\GL}\) truly corresponds to \(\GL\), in other words we have the following theorem.\footnote{In fact \(\g{\GL}\) is an inessential variation of the sequent calculus presented at \cite{Valentini}.}

\begin{theorem}[\cite{Valentini}]
  \label{calculus-GL-sound-and-complete}
  For any sequent \(\Gamma \Rightarrow \Delta\) we have that \(\GL \vdash \bigwedge \Gamma \to \bigvee \Delta\) if and only if \(\g{\GL} \vdash \Gamma \Rightarrow \Delta\).
\end{theorem}

\begin{figure}
  \[
    \AxiomC{\(\Gamma \Rightarrow \Delta\)}
    \RightLabel{\(\wk\)}
    \UnaryInfC{\(\Gamma, \Gamma' \Rightarrow \Delta, \Delta'\)}
    \DisplayProof
    \quad
    \AxiomC{\(\phi, \phi, \Gamma \Rightarrow \Delta\)}
    \RightLabel{\(\ctr\)}
    \UnaryInfC{\(\phi, \Gamma \Rightarrow \Delta\)}
    \DisplayProof
    \quad
    \AxiomC{\(\Gamma \Rightarrow \phi, \phi, \Delta\)}
    \RightLabel{\(\ctr\)}
    \UnaryInfC{\(\Gamma \Rightarrow \phi, \Delta\)}
    \DisplayProof
  \]
  \caption{\(\GL\) structural rules}
  \label{fig:GL-structural-rules}
  \separator
\end{figure}

Thanks to Theorem~\ref{calculus-GL-sound-and-complete} we know that proofs in \(\g{\GL}\) are witnesses of validities of formulas in \(\GL\).
However, a quick observation to Solovay's original proof will show that Kripke models are used in the opposite way: as witnesses of the non-validity of formulas (i.e., as countermodels).
For that reason we need to define refutations, which are the proof theoretical equivalent of countermodels.
However, first we remember that \(\g{\GL}\) is well-behaved with respect to structural\footnote{Note that the notion of structural in proof theory is distinct from the notion of structural property of Kripke models that we discussed in the introduction.} rules (see Figure~\ref{fig:GL-structural-rules}), as usual the lemma is proved by induction on the height of the proof.

\begin{lemma}\label{lem:Wk-Ctr-Invert}
  We have the following:
  \begin{enumerate}
    \item \((\wk)\) is admissible preserving height in \(\g{\GL}\),
    \item \((\botR)\), \((\toL)\) and \((\toR)\) are invertible preserving height in \(\g{\GL}\),
    \item \((\ctr)\) is admissible preserving height in \(\g{\GL}\).
  \end{enumerate}
\end{lemma}

\begin{figure}
  \[
    \AxiomC{\(\phi \to \psi, \Gamma \Rightarrow \phi,\Delta\)}
    \RightLabel{\(\toLsat_0\)}
    \UnaryInfC{\(\phi \to \psi, \Gamma \Rightarrow \Delta\)}
    \DisplayProof
    \qquad
    \AxiomC{\(\psi,\phi \to \psi, \Gamma \Rightarrow \Delta\)}
    \RightLabel{\(\toLsat_1\)}
    \UnaryInfC{\(\phi \to \psi, \Gamma \Rightarrow \Delta\)}
    \DisplayProof
    \qquad
    \AxiomC{\(\phi, \Gamma \Rightarrow \psi, \phi \to \psi,\Delta\)}
    \RightLabel{\(\toRsat\)}
    \UnaryInfC{\(\Gamma \Rightarrow \phi \to \psi,\Delta\)}
    \DisplayProof
  \]
  \caption{Rules for proof search}
  \label{fig:rules-ps-GL}
\end{figure}

To define the proof search we will need the notion of saturation below and the rules of Figure~\ref{fig:rules-ps-GL}.

\begin{definition} \label{def:saturation-GL}
  Given a sequent \(S = (\Gamma \Rightarrow \Delta)\) we say that
  \begin{itemize}
    \item \(\phi\) is \emph{left-saturated in \(S\)} if
    \begin{enumerate}
      \item \(\phi\) is a propositional variable or a \(\nec\)-formula, or
      \item \(\phi = \phi_0 \to \phi_1\) and \(\phi_0 \in \Delta\) or \(\phi_1 \in \Gamma\).
    \end{enumerate}
    \item \(\phi\) is \emph{right-saturated in \(S\)} if
    \begin{enumerate}
      \item \(\phi\) is a propositional variable or a \(\nec\)-formula, or
      \item \(\phi = \phi_0 \to \phi_1\) and \(\phi_0 \in \Gamma\) and \(\phi_1 \in \Delta\).
    \end{enumerate}
    \item \(S\) is said to be \emph{saturated} if every \(\phi \in \Gamma\) is left-saturated in \(S\) and every \(\psi \in \Delta\) is right-saturated in \(S\).
  \end{itemize}
  We define the saturation complexity of a sequent \(S = (\Gamma \Rightarrow \Delta)\) as
  \[
    \satc{S} = \sum \set{|\phi| \mid (\phi \in \Gamma \text{ and } \phi \text{ not left-saturated in }S) \text{ or }(\phi \in \Delta \text{ and } \phi \text{ not right-saturated in }S)},
  \]
  where \(|\phi|\) is the complexity of the formula.
  In words, the saturation complexity of a sequent is the sum of the complexities of the formulas in \(S\) which are not saturated.
\end{definition}

\begin{definition}[Proof Search Game]
  Given a sequent \(S\), we define the rule set for proof search in \(\GL\) denoted \(\ps_{\GL}\) as follows.
  \begin{enumerate}
    \item The positions are sequents \(S\).
    \item We say that \(S'\) is a move at \(S\) if
      \begin{enumerate}
        \item \(S\) is not an extended axiomatic sequent, and
        \item \((S',S)\) is an instance of one of the rules \((\toLsat_0), (\toLsat_1), (\toRsat)\) with a non-saturated principal formula, or \(S\) is saturated and \((S',S)\) is a maximal instance of \((\modal{\GL})\).
      \end{enumerate}
    \item The positions of player \(\pI\) are the sequents which are saturated and not extended axiomatic.
      Notice that then, positions of player \(\pII\) are the sequents which are either non-saturated or extended axiomatic.
  \end{enumerate}
  It will be common to call player \(\pI\) \emph{verifier} and player \(\pII\) \emph{refuter}.
  A winning strategy for refuter at game \(\ps_{\GL} @ S\) is also called a \emph{refutation of \(S\)}.
\end{definition}
We are going to show that for any sequent \(S\) the game \(\ps_{\GL} @ S\) is a short game.
Define the following measures on sequents:
\begin{align*}
  &|\Gamma \Rightarrow \Delta|_{\nec} = |\sub(\Gamma \Rightarrow \Delta)\setminus \set{\nec \phi \mid \nec \phi \in \Gamma}| \text{, i.e.,\ the cardinality of the set }\sub(\Gamma \Rightarrow \Delta)\setminus \set{\nec \phi \mid \nec \phi \in \Gamma}\text{;}\\
  &|S|_{\GL} = \omega |S|_{\nec} + \satc{S}.
\end{align*}
We have the following.
\begin{lemma}\label{GL-property-of-the-measure}
  If \(S,S'\) are sequents such that \(S'\) is a move at \(S\) then \(|S|_{\GL} > |S'|_{\GL}\).
\end{lemma}
\begin{proof}
  Since \(S'\) is a move at \(S\) we have that \(S\) is not an extended axiomatic sequent and that \((S',S)\) is an instance of at least one of \((\toLsat_0)\), \((\toLsat_1)\), \((\toRsat)\) or \((\modal{\GL})\).
  If it is an instance of one of the first three rules, we have that the principal formula is not saturated.
  It is clear from the shape of the rules \((\toLsat_0)\), \((\toLsat_1)\) and \((\toRsat)\) that then \(|S|_{\nec} \geq |S'|_{\nec}\), as no \(\nec\)-formulas are eliminated from the left side of the sequent and the set of subformulas of both sequents is the same.
  Since the principal formula was not saturated, but it becomes saturated after the application of the rule (possibly with the introduction of non-saturated formulas whose sum of complexities is smaller than the complexity of the original formula), we have that \(\satc{S} > \satc{S'}\).
  Both facts together give the desired \(|S|_{\GL} > |S'|_{\GL}\).

  Finally, assume that \((S',S)\) is an instance of \((\modal{\GL})\).
  Then the rule instance is of shape
  \[
    \AxiomC{\(\necd \Sigma, \nec \phi \Rightarrow \phi\)}
    \RightLabel{\(\modal{\GL}\)}
    \UnaryInfC{\(\nec \Sigma, \Gamma \Rightarrow \nec \phi, \Delta\)}
    \DisplayProof
  \]
  where \(S'\) is the premise and \(S\) the conclusion.
  As the rule application needs to be maximal we have that \(\Gamma\) does not contain \(\nec\)-formulas and since \(S\) is not an extended axiomatic sequent we have that \(\nec \phi \not \in \nec \Sigma\).
  This implies that \(|S|_{\nec} > |S'|_{\nec}\), so \(|S|_{\GL} > |S'|_{\GL}\).
\end{proof}

Using the measure \(|\cdot|_{\GL}\), it is then easy to establish that \(\ps_{\GL} @ S\) is a short game.

\begin{theorem}\label{short-game-psGL}
  For any sequent \(S\), \(\ps_{\GL}@ S\) is a short game.
  As a corollary, either verifier or refuter have a winning strategy in \(\ps_{\GL} @ S\).
\end{theorem}
\begin{proof}
  Suppose that \((S_i)_{i < \omega}\) is an infinite run of \(\ps_{\GL} @ S\).
  Then \(|S_0|_{\GL} > |S_1|_{\GL} > \cdots > |S_n|_{\GL} > \cdots\) is an infinite strictly descending chain of ordinals, which is impossible.
  So \(\ps_{\GL} @ S\) has no infinite runs.

  To show that for any \(S\),  \(\ps_{\GL} @ S\) has finitely many subpositions, we proceed by induction on \(|S|_{\GL}\).
  So assume it holds for any \(S'\) such that \(|S'|_{\GL} < |S|_{\GL}\).
  Notice that the set \(\set{S' \mid S' \text{ is a move at }S}\) is finite: if \(S\) is non-saturated there is a move for each non-saturated formula and if it is saturated there is a move for each distinct \(\nec\)-formula of the right of \(S\).
  Then let \(S_0,\ldots,S_n\) be the moves at \(S\), by the previous lemma we have that \(|S_i|_{\GL} < |S|_{\GL}\).
  By the induction hypothesis \(S_i\) has also finitely many subpositions for \(i < n\), and then it is clear that the \(S\) has also finitely many subpositions.
\end{proof}

We are ready to show that from either player's winning strategy we can build either a proof or a refutation of the sequent.

\begin{theorem}
  \label{fundamental-theorem-of-refutations-GL}
  For any sequent \(S\) we have that
  \begin{enumerate}
    \item \(\g{\GL} \vdash S\) if and only if verifier has a winning strategy in \(\ps_{\GL} @ S\).
    \item \(\g{\GL} \not\vdash S\) if and only if refuter has a winning strategy in \(\ps_{\GL} @ S\), i.e., if \(S\) has a refutation.
  \end{enumerate}
\end{theorem}
\begin{proof}
  Proof of 1. For the left to right direction can just use induction on the height of the proof and height preserving admissibility of weakening.
  For the right to left direction it suffices to use induction on the conversely wellfounded order \(\leadsto\) and admissibility of contraction.
  
  Proof of 2. This follows immediately from the first part.
\end{proof}

\subsection{Arithmetical completeness for \(\GL\)}
\label{sec:arithmetical-completeness-GL}

For the remainder of Subsection~\ref{sec:arithmetical-completeness-GL} let us fix an arithmetical theory \(T\) and a predicate \(\prv\) fulfilling the conditions stated at Theorem~\ref{arithmetical-completeness-GL}.

Now, assume we have a formula \(\phi\) such that \(\GL \not\vdash \phi\), i.e., \(\g{\GL} \not \vdash { \Rightarrow \phi}\).
This implies that there is a refutation of \(S\), let us denote this winning strategy for refuter by \(s_{\pII}\).
Since \(\ps_{\GL} @ S\) is short we have that \(\set{w \in \run(S| s_{\pII}) \mid \last{w} \in P_{\pI}}\) of \(s_{\pII}\)-compatible runs whose last sequent is a position of verifier is also finite, let it be \(\set{w_1, \ldots, w_N}\).
Additionally, assume that we have chosen the indices such that if the run \(w_y\) extends the run \(w_x\) then \(x < y\) (this can be done easily because the set of runs with \(\leadsto^+\) is a finite tree).
Define \(W = \set{0,\ldots,N} \subseteq \mathbb{N}\), and define the relation \(x \sqsubset y\) for \(x,y \in W\) as
\begin{enumerate}
  \item \(x = 0 \) and \(y = 1\), or
  \item \(x,y > 0\), \(w_x = (S_0,\ldots,S_n)\), \(w_y = (S_0,\ldots,S_n,S_{n+1}, \ldots,S_m)\) for \(n < m\) and there is no position of verifier in \(S_{n+1}, \ldots,S_{m-1}\).
\end{enumerate}
For \(x \in W \setminus \set{0}\) let us denote the last sequent of \(w_x\) by \(S_x = \Gamma_x \Rightarrow \Delta_x\) and let us define \(S_0 := \Gamma_0 \Rightarrow \Delta_0 := {\Rightarrow \phi}\).
Since if \(w_y\) extends \(w_x\) we have that \(x < y\) we can see that the final position of \(w_1\) is the only position belonging to verifier.
This means that \(S_1\) is obtained from \(S_0\) by applying the rules \((\toL_0^{\mathrm{sat}})\), \((\toL_1^{\mathrm{sat}})\) and \((\toR^{\mathrm{sat}})\) so \(\Delta_0 \subseteq \Delta_1\).
This property will be useful in what follows.

\label{remarks-about-GL-tree}
Intuitively, the points in the structure $(W, {\sqsubset})$ can be viewed as finite branches in the proof search game which end in a saturated sequent.
  If \(x \sqsubset y\) then \(w_y\) is, in a sense, an \emph{immediate successor} of \(w_x\): from \(w_x\) to \(w_y\) verifier only does one move and the rest of the moves by refuter are unique since they are directed by \(s_{\pII}\).
In fact it is easy to check that \((W,\sqsubset^+)\) is a finite rooted tree with root \(0\) and \(\sqsubset\) is the corresponding immediate successor relation.
For the semantic-minded reader, we want to notice that this structure forms a finite tree akin to the countermodel for $\GL$ that is needed in Solovay's completeness proof.
However, note that in our case the tree is not transitive; instead, the ``transitivity'' is a consequence of the rules of the sequent calculus.
In fact, the corresponding Kripke model to this structure would not even be a model of \(\GL\), it will only satisfy enough instances of the axioms of \(\GL\) to guarantee the unprovability of \(\phi\) in \(\GL\).
Since our Solovay function will move through the \(\sqsubset\)-relation and not through the \(\sqsubset^+\)-relation we will embed the non-transitive tree in arithmetic, and not its transitive closure.

The following lemma shows that, even if the tree does not have the usual structural properties, it carries all the necessary information for the refutation of \(\phi\) in \(\GL\) by virtue of the shape of the rules of \(\g{\GL}\).
\begin{lemma}\label{properties-of-refutations-GL}
  We have the following for \(x \in W \setminus \set{0}\),
  \begin{enumerate}
    \item \(\bot \not \in \Gamma_x\), \(\Gamma_x \cap \Delta_x = \varnothing\) and \(\Gamma_x \Rightarrow \Delta_x\) is saturated,
    \item \(\nec \psi \in \Delta_x\) implies that there is an \(y\) such that \(x \sqsubset y\) and \(\psi \in \Delta_y\).
    \item \(\nec \psi \in \Gamma_x\) and \(x \sqsubset^+ y\) implies \(\nec \psi, \psi \in \Gamma_y\),
  \end{enumerate}
\end{lemma}
\begin{proof}
  The first property is thanks to the positions of verifier being saturated sequents which are not extended axiomatic.

  The second property follows from the fact that if \(\nec\psi \in \Delta_x\) then \(S_x = (\nec \Sigma, \Gamma \Rightarrow \nec \psi, \Delta)\) for some multisets \(\Sigma, \Gamma, \Delta\) such that \(\Gamma\) does not contain \(\nec\)-formulas.  Then, as \(\hat{w}_x = S_x\) is a position of verifier, it follows that \(S_x\) is saturated and non-extended axiomatic.
  This means that if we define \(S' := (\necd \Gamma, \nec \psi \Rightarrow \psi)\) we will have that \(S \leadsto S'\), so adding \(S'\) to the end of \(w_x\) will still be a \(s_{\pII}\)-consistent run of the game.
  Following the strategy \(s_{\pII}\), we can extend this run until we obtain a position of verifier (since \(s_{\pII}\) is a winning strategy for refuter we cannot get stuck, and since the game is short we must eventually arrive to a position belonging to verifier).
  This extended run will have the form \(w_y\) for some \(y \leq N\) and we will have that \(x \sqsubset y\).
  Also, since the moves of refuter only add formulas without deleting any we will have that \(\psi \in \Delta_y\).

  For the third property it suffices to show that \(\nec \psi \in \Gamma_x\) and \(x \sqsubset y\) implies \(\nec \psi, \psi \in \Gamma_y\).
  This is a direct consequence of the moves of refuter only adding formulas to the sequent and the shape of the move of verifier, which are maximal instances of \(\modal{\GL}\).
  
\end{proof}

We will use these properties of refutations throughout the remainder of Subsection~\ref{sec:arithmetical-completeness-GL}, sometimes without explicitly mentioning them.
With these preliminaries in place, we are ready to prove the desired result.

\subsubsection{The heart of the proof}
In this subsubsection, we show how to construct the arithmetical interpretation \((\cdot)^*\) such that \(T \not \vdash \phi^*\).
For this we will assume the existence of some arithmetical sentences, called \emph{Solovay sentences}, \((\lambda_0,\ldots,\lambda_N)\) such that
\begin{multicols}{2}
\begin{enumerate}\label{solovay-properties-GL}
  \item (\(\mathrm{S1}\)) \(T \vdash \bigvee_{x \leq N} \lambda_x\),
  \item (\(\mathrm{S2}\)) \(T \vdash \neg(\lambda_x \wedge \lambda_y)\) for \(x \neq y\),
  \item (\(\mathrm{S3}\)) \(T \vdash \lambda_x \to \posg \lambda_y\) for \(x \sqsubset y\),
  \item (\(\mathrm{S4}\)) \(T \vdash \lambda_x \to \necg\bigvee_{y \mid x \sqsubseteq^+ y} \lambda_y\),
  \item (\(\mathrm{S5}\)) \(T \vdash \lambda_x \to \necg\neg \lambda_x\) for \(x \neq 0\).
\end{enumerate}
\end{multicols}
In the next subsubsection, we will show how to construct such sentences using the ideas of Solovay.
The following is an immediate consequence of the properties of Solovay sentences.
\begin{lemma}
  We have that \(\mathbb{N} \vDash \lambda_0\).
\end{lemma}
\begin{proof}
  Since \(T \vdash \bigvee_{x \in W} \lambda_x\) and \(T\) is sound we have that \(\mathbb{N} \vDash \bigvee_{x \in W} \lambda_x\), so \(\mathbb{N} \vDash \lambda_x\) for some \(x \in W\), assume that \(x \neq 0\).
  Then \(\mathbb{N} \vDash \lambda_x \to \necg \neg \lambda_x\) by (S5) and the soundness of \(T\), which implies \(\mathbb{N} \vDash \necg \neg \lambda_x\).
  The soundness of \(\prv\) gives that \(T \vdash \neg \lambda_x\) and soundness of \(T\) allow us to conclude \(\mathbb{N} \not\vDash \lambda_x\), a contradiction.
\end{proof}

We define the arithmetical interpretation
\[
  p^* := \bigvee_{x \mid p \in \Gamma_x} \lambda_x,
\]
which will be the desired interpretation such that \(T \nvdash \phi^*\).
To show this, we first need to prove a truth lemma.

\begin{lemma}[Truth Lemma]\label{truth-lemma-GL}
  For any \(x \in W \setminus \set{0}\), we have that
  \begin{align*}
    &\psi \in \Gamma_x \text{ implies } T \vdash \lambda_x \to \psi^*, 
    &&\psi \in \Delta_x \text{ implies } T \vdash \lambda_x \to \neg \psi^*.
  \end{align*}
\end{lemma}
\begin{proof}
  We proceed by induction on the complexity of \(\phi\).
  If \(\psi = p\) for some propositional variable \(p\), the first result is a trivial consequence of logical reasoning and the second result follows from \((\mathrm{S2})\) straightforwardly.

  If \(\psi = \bot\) then we have that \(\bot \not \in \Gamma_x\), so the first result is vacously satisfied and the second result becomes \(T \vdash \lambda_x \to \neg \bot\), which holds straightforwardly.

  If \(\phi = \phi_0 \to \phi_1\) and \(\phi_0 \to \phi_1 \in \Gamma_x\) we have by saturation that either \(\phi_0 \in \Delta_x\) or \(\phi_1 \in \Gamma_x\).
  In the first case the induction hypothesis gives \(T \vdash \lambda_x \to \neg \phi^*_0\) while the second gives \(T \vdash \lambda_x \to \phi^*_1\).
  In both cases, \(T \vdash \lambda_x \to (\phi^*_0 \to \phi^*_1)\) by logical reasoning.
  Now assume that \(\phi_0 \to \phi_1 \in \Delta_x\), by saturation  \(\phi_0 \in \Gamma_x\) and \(\phi_1 \in \Delta_x\), so by induction hypothesis \(T \vdash \lambda_x \to \phi^*_0\) and \(T \vdash \lambda_x \to \neg \phi_1^*\).
  By logical reasoning \(T \vdash \lambda_x \to \neg (\phi^*_0 \to \phi^*_1)\), as desired.

  Case \(\phi = \nec \phi_0\), assume that \(\nec \phi_0 \in \Gamma_w\).
  We have the following reasoning.
  \begin{align*}
    &T \vdash \lambda_y \to \phi^*_0 \text{ for \(x \sqsubset^+ y\)},
    &&\text{using the I.H.\ since \(x \sqsubset^+ y\) and \(\nec \phi_0 \in \Gamma_x\) implies \(\phi_0 \in \Gamma_y\);} \\
    &T \vdash \bigvee_{y \mid x \sqsubset^+ y} \lambda_y \to \phi^*_0, 
    &&\text{by logical reasoning;} \\
    &T \vdash \necg\left(\bigvee_{y \mid x \sqsubset^+ y } \lambda_y\right) \to \necg\phi^*_0,
    &&\text{by \((\NEC[\necg])\) and \((\Kax[\necg])\)}; \\
    &T \vdash \necg\left(\bigvee_{y \mid x \sqsubseteq^+ y } \lambda_y\right) \to \necg\phi^*_0,
    &&\text{by (\(\mathrm{S}5\)) since \(x \neq 0\)},\\
    &T \vdash \lambda_x \to (\nec\phi_0)^* 
    &&\text{by (\(\mathrm{S}4\))}.
  \end{align*}
  Finally, assume that \(\nec \phi_0 \in \Delta_x\).
  We know that there is a \(y\) such that \(x \sqsubset y\) and \(\phi_0 \in \Delta_y\).
  We have the following reasoning:
  \begin{align*}
    &T \vdash \lambda_y \to \neg \phi_0^*
    &&\text{by the I.H.;} \\
    &T \vdash \phi_0^* \to \neg \lambda_y
    &&\text{by logical reasoning;} \\
    &T \vdash \necg\phi_0^* \to \necg\neg \lambda_y
    &&\text{by \((\NEC[\necg])\) and \((\Kax[\necg])\)}; \\
    &T \vdash \posg\lambda_y \to \neg\necg\phi_0^*
    &&\text{by logical reasoning;} \\
    &T \vdash \lambda_x \to \neg (\nec \phi_0)^*
    &&\text{using \((\mathrm{S3})\) since \(x \sqsubset y\).}
  \end{align*}
  This establishes the desired result.
\end{proof}


We now use the Truth Lemma to establish arithmetical completeness (Theorem~\ref{arithmetical-completeness-GL}).
Recall that \(\set{\phi} = \Delta_0 \subseteq \Delta_1\) and \(0 \sqsubset 1\).
Assume that \(T \vdash \phi^*\) and we get a contradiction as follows.
\begin{align*}
  &T \vdash \neg \lambda_1,
  &&\text{since }T \vdash \phi^* \text{ and by the Truth Lemma }T \vdash \lambda_1 \to \neg \phi^*; \\
  &T \vdash \necg\neg \lambda_1,
  &&\text{by \((\NEC[\necg])\)}; \\
  &T \vdash \neg \lambda_0,
  &&\text{by \((\mathrm{S}3)\)}; \\
  &\mathbb{N} \not\vDash  \lambda_0,
  &&\text{by soundness of \(T\)};
\end{align*}
thus we have obtained the desired contradiction since \(\mathbb{N} \vDash \lambda_0\).

\begin{example}\label{example-GL}
  Let \(T\) be a sound arithmetical theory with a \(\Sigma_1\)-formula enumerating its axioms.
  Then it is possible to construct a sound Hilbert-Bernays provability predicate \(\prv[T](x)\) that internalizes naturally the notion of having a Hilbert proof from the axioms of \(T\)~\cite{Feferman1960}.
  Examples of theories for which this applies are \(\PA\), \(\EA\), the theories \(\IS{n}\) for \(n \geq 1\) among others.
  Then if \(T\) is any of these theories we can conclude that \(\GL\) is arithmetically complete with respect to \((\cdot)^*_{\prv[T]}\) in \(T\).
\end{example}

\begin{remark}[Known generalizations]
  Let us quickly mention some generalizations which are known in the literature.
  \begin{enumerate}
    \item It is possible to change the soundness of \(T\) for \(\Sigma_1\)-soundness and the soundness of \(\prv\) for the property that for any \(n < \omega\), \(T \not\vdash \necg^n \bot\), where \(\necg^ 0 \bot := \bot\) and \(\necg^{n+1} \bot := \necg \necg^n \bot\).
    \item The previous generalization implies that it suffices that \(T\) is \(\Sigma_1\)-sound and \(\prv\) is sound.
    \item It is possible to change the base theory and consider single-conclusion Hilbert-Bernays provability predicate for \(T\) in \(B\), if \(B\) fulfils some conditions such as being at least as strong as \(\EA\).
    \item If \(\prv(x) = \exists p(\prf(p, x))\) is a sound Hilbert-Bernays provability predicate for \(T\) but \(\prf(p,x)\) is not single conclusion we could define \(\prf^1(p,x) := p = \overline{2}^q \times \overline{3}^x \wedge \prf(q,x)\), it is straightforward to check that \(T\) proves \(\prf^1\) to be single-conclusion.
    Define \(\prv^1(x)\) as the strict \(\Sigma_1\)-formula \(\exists p (\prf^1(p, x)) \), then it is easy to check that \(T \vdash \forall x (\prv(x) \leftrightarrow \prv^1(x))\), from which it is easy to see that \(T \vdash (\phi)^*_{\prv} \leftrightarrow (\phi)^*_{\prv^1}\) for any arithmetical interpretation \((\cdot)^*\).
    Then \(\prv^1\) is a sound single-conclusion Hilbert-Bernays provability predicate, and then we can conclude that the arithmetical completeness also follows for provability predicates whose proof predicate is not single conclusion.
  \end{enumerate}
\end{remark}

\subsubsection{Construction of the Solovay function}

We are going to define the sentences \(\lambda_0, \ldots, \lambda_N\) fulfilling the properties \((\mathrm{S}1)\), \ldots, \((\mathrm{S}5)\) of page~\pageref{solovay-properties-GL}.
The intuition behind Solovay's construction is as follows.
Imagine that, inside \(T\), we define a function \(f\) whose range is \(W\).
It starts at \(0\) and if \(f(n) = x\), \(x \sqsubset y\) and \(n\) is the G\"odel code of a proof in \(\prv\) that \(f\) does not stay at \(y\) then \(f(n+1) = y\).
Since this function is climbing through a finite tree it is easy to see that it has a limit, then \(\lambda_x\) claims that the limit of the Solovay function is \(x\).

To avoid the complications of defining a self-referential function we will follow the ideas of~\cite{japaridze1986modal} and of~\cite{dejongh1991proof} and use the uniform diagonal lemma for \(T\).
 With this purpose in mind we will define auxiliary sentences \(\rho_x\) for \(x \in W\) whose meaning can be understood as claiming that the Solovay function passes through \(x\).
We define the arithmetical formulas \(\rho_x(p_0,\ldots,p_N)\) for each $x \in W$ with free variables $p_0,\ldots,p_N$ as follows
\begin{align*}
  &\rho_0(p_0,\ldots,p_N) := \top, && \\
  &\rho_x(p_0,\ldots,p_N) := \rho_y(p_0,\ldots,p_N) \wedge \prv(\dot{\neg} p_x) \wedge \bigwedge_{z \mid y \sqsubset z \neq x} \prv(\dot\neg p_x) < \prv(\dot\neg p_z)
  &&\text{where }y \sqsubset x.
\end{align*}
In the definitions above we define \(\prv(\dot{\neg} x)\) as a strict \(\Sigma_1\)-formula equivalent in \(T\) to the formula \(\exists p, y (\mathrm{Neg}(x,y) \wedge \prf(p,y))\), i.e., as saying that the negation of \(x\) is provable.
Note that such a strict \(\Sigma_1\)-formula exists since \(\mathrm{Neg}(x,y)\) is \(\Delta_0\) and then the formula above is \(\Sigma_1\).

The \(\rho_x\)'x are well-defined since each \(x \in W \setminus \set{0}\) has a unique \(y \in W\) such that \(y \sqsubset x\) as \(\sqsubset\) is the immediate successor relation of a tree.\footnote{Note that this is the exact reason why we decide to keep the tree structure of the tree of runs of the proof search. It would be possible to instead of using the runs of the proof search, use just the positions and instead of having exactly one \(y\) we could have multiple of them for each \(x > 0\). This, even if possible, makes the following proofs harder so for this reason we will stick to the tree.
}
We consider the system of formulas
\[
  \rho_x(p_0,\ldots,p_N) \wedge \bigwedge_{y \mid x \sqsubset y} \con(p_y)\qquad \text{for }x \in W,
\]
where we define \(\con(x) := \neg \prv(\dot{\neg} x)\), i.e.,\ the consistency predicate of \(\prv\).
Given an arithmetical sentence \(\alpha\) we will write \(\posg \alpha\) to mean \(\con(\Gdl{\alpha})\).

Thanks to the uniform diagonal lemma in \(T\), we obtain arithmetical sentences  \(\lambda_0,\ldots,\lambda_N\) such that
\[
  T \vdash \lambda_x \leftrightarrow \rho_x \wedge \bigwedge_{y \mid x \sqsubset y} \posg \lambda_y, \qquad \text{for }x \in \set{0,\ldots,N},
\]
where \(\rho_x\) is just \(\rho_x := \rho_x(\ulcorner \lambda_0\urcorner, \ldots,\ulcorner\lambda_N\urcorner)\) for \(x \in W\).
Additionally, by definition of the \(\rho_x\)'s we have
\[
  T \vdash \rho_x \leftrightarrow \rho_y \wedge \necg \neg \lambda_x \wedge \bigwedge_{z \mid y \sqsubset z \neq x} \necg \neg \lambda_x < \necg \neg \lambda_z \qquad \text{where }y \sqsubset x.
\]
We will also assume that for \(x \neq y\) we have \(\lambda_x \neq \lambda_y\), i.e.,\ that for different elements of \(W\) the Solovay sentences are not syntactically equal.
Otherwise, if \(\lambda_x = \lambda_y\) for \(x \neq y\), we could redefine \(\lambda_y\) as \(\lambda_x \wedge \top\) and as there are only finitely many Solovay sentences this modification suffices.
As straightforward consequences of the definitions we have the following propositions.

\begin{proposition}[S3]
  We have that \(T \vdash \lambda_x \to \posg \lambda_y\) for \(x \sqsubset y\).
\end{proposition}

\begin{proposition}[S5]
  We have that \(T \vdash \lambda_x \to \necg \neg \lambda_x\) for \(x \in W \setminus \set{0}\).
\end{proposition}

Now let us move to the proof of (S2).
We start by showing that if \(x\) and \(y\) do not belong to the same branch then the Solovay function cannot pass through both of them.

\begin{lemma}
  \label{incomparable-absurd-GL}
  If \(x \not\sqsubseteq^+ y\) and \(y \not\sqsubseteq^+ x\) then \(T \vdash \neg(\rho_x \wedge \rho_y)\).
\end{lemma}
\begin{proof}
  Since \(x \not\sqsubseteq^+ y\) and \(y \not\sqsubseteq^+ x\) we have that \(x,y\) must belong to different branches of the tree \((W,\sqsubset^+)\).
  Then, there must be \(z_0,z_1,z_2\) such that \(z_0 \sqsubset z_1 \sqsubseteq^+ x\) and \(z_0 \sqsubset z_2 \sqsubseteq^+ y\) and \(z_1 \neq z_2\).
  By definition we can see that \(T \vdash \rho_x \to \rho_{z_1}\)  and \(T \vdash \rho_y \to \rho_{z_2}\), so it suffices to show that \(T \vdash \neg(\rho_{z_1} \wedge \rho_{z_2})\).
  Since \(z_0 \sqsubset z_1 \neq z_2\) and \(z_0 \sqsubset z_2 \neq z_1\) then \(T \vdash \rho_{z_1} \to \necg \neg \lambda_{z_1} < \necg \neg \lambda_{z_2}\) and \(T \vdash \rho_{z_2} \to \necg \neg \lambda_{z_2} < \necg \neg \lambda_{z_1}\).
  Since, using that \(T\) shows linearity of \(\preceq\), \(T \vdash \neg( \necg \neg \lambda_{z_1} < \necg \neg \lambda_{z_2} \wedge \necg \neg \lambda_{z_2} < \necg \neg \lambda_{z_1})\) we can conclude the desired result.
\end{proof}

Also, if \(y\) is a node strictly above \(x\) then it is not the case that the limit of the Solovay function is at \(x\) and passes through \(y\).

\begin{lemma}
  \label{limit-absurd-GL}
  If \(x \sqsubset y\) then \(T \vdash \neg(\lambda_x \wedge \rho_y)\).
\end{lemma}
\begin{proof}
  We have that \(T \vdash \lambda_x \to \posg \lambda_y\) and \(T \vdash \rho_y \to \necg \neg \lambda_y\), from which the desired result follows.
\end{proof}

We can conclude that the Solovay function as at most one limit.

\begin{theorem}[S2]\label{property-S2-GL}
  We have that \(T \vdash \neg(\lambda_x \wedge \lambda_y)\) for \(x \neq y\).
\end{theorem}
\begin{proof}
  Given \(x \neq y\) we have that either \(x \not\sqsubseteq^+ y\) and \(y \not\sqsubseteq^+ x\), or \(x \sqsubset^+ y\), or \(y \sqsubset^+ x\).
  If \(x \not\sqsubseteq^+ y\) and \(y \not\sqsubseteq^+ x\) we obtain the desired result since \(T \vdash \lambda_x \to \rho_x\), \(T \vdash \lambda_y \to \rho_y\) and \(T \vdash \neg(\rho_x \wedge \rho_y)\) by Lemma~\ref{incomparable-absurd-GL}.
  We treat the case when \(x \sqsubset^+ y\), if \(y \sqsubset^+ x\) the reasoning is analogous.
  Note that if \(x \sqsubset^+ y\) then there is a \(z\) such that \(x \sqsubset z \sqsubseteq^+ y\).
  By definition of the \(\rho\)s, it is easy to see that \(T \vdash \rho_y \to \rho_z\).
  Since \(T \vdash \lambda_y \to \rho_y\) by definition of the \(\lambda\)s and \(T \vdash \neg(\lambda_x \wedge \rho_z)\) by Lemma~\ref{limit-absurd-GL} we obtain the desired result.
\end{proof}

The final property we need to show is that if the Solovay function passes through a node \(x\) then the limit must be either \(x\) or a node above \(x\), from this all the rest of Solovay conditions can be derived.
To show this first we need the following lemma.

\begin{lemma}
  \label{finding-minimum-GL}
  We have that \(T \vdash \left(\bigvee_{y \mid x \sqsubset y} \necg \neg \lambda_y\right) \to \left(\bigvee_{y \mid x \sqsubset y} \necg \neg \lambda_y \wedge \bigwedge_{z \mid x \sqsubset z \neq y} \necg \neg \lambda_y < \necg \neg \lambda_z\right)\).
\end{lemma}
\begin{proof}
  We are going to show that for any finite set \(\Gamma\) of (different) sentences we have that for \(\gamma \in \Gamma\)
  \[
    T \vdash \necg \gamma \to \left(\bigvee_{\beta \in \Gamma} \necg \beta \wedge \bigwedge_{\beta' \in \Gamma \setminus \set{\beta}} \necg \beta < \necg \beta'\right).
  \]
  Note that this can be understood as saying that, inside \(T\), if one of the sentences in \(\Gamma\) is \(\prv\)-provable then there is a unique sentence in \(\Gamma\) which has a minimum \(\prv\)-proof.
  This is enough to establish the desired lemma since we know that for \(x \neq y\), \(\lambda_x \neq \lambda_y\).
  
  By induction on the size of \(\Gamma\), if \(\Gamma\) is empty then we have \(T \vdash \bot \to \bot\), as desired.
  Assume it holds for sets of size \(n\) and \(\Gamma\) is of size \(n + 1\), let \(\gamma \in \Gamma\) and define \(\Gamma' = \Gamma \setminus \set{\gamma}\).
  We know that, using \(\Sigma_0\)-minimization in \(T\) and that \(\prf\) is single conclusion, for any \(\alpha \in \Gamma'\) we have \(T \vdash \necg \gamma \to \necg \alpha < \necg \gamma \vee \necg \gamma < \necg \alpha\).
  Then, we have that
  \[T \vdash \necg \gamma \to \left(\necg \gamma \wedge \bigwedge_{\alpha \in \Gamma'} \necg \gamma < \necg \alpha\right) \vee\bigvee_{\alpha \in \Gamma'} \necg \alpha < \necg \gamma.\]
  Note that, thanks to this, it suffices to show the following assertions:
  \begin{align*}
    \tag{i}
    &T \vdash \left(\necg \gamma \wedge \bigwedge_{\alpha \in \Gamma'} \necg \gamma < \necg \alpha\right) \to \left(\bigvee_{\beta \in \Gamma} \necg \beta\wedge \bigwedge_{\beta' \in \Gamma \setminus \set{\beta}} \necg \beta < \necg \beta'\right) \\
    \tag{ii}
    &T \vdash \necg \alpha < \necg \gamma \to \left(\bigvee_{\beta \in \Gamma} \necg \beta\wedge \bigwedge_{\beta' \in \Gamma \setminus \set{\beta}} \necg \beta < \necg \beta'\right) \quad \text{ for }\alpha \in \Gamma'.
  \end{align*}

  (i) is clearly provable by propositional reasoning as the antecedent is a disjunct of the consequent, so we focus in proving (ii).
  Let \(\alpha \in \Gamma'\) we have that \(T \vdash \necg \alpha < \necg \gamma \to \necg \alpha\).
  Using the induction hypothesis on \(\alpha \in \Gamma'\) we obtain that 
  \[T \vdash \necg \alpha < \necg \gamma \to \left(\bigvee_{\beta \in \Gamma'} \necg \beta \wedge \left(\bigwedge_{\beta \in \Gamma' \setminus \set{\beta}} \necg \beta < \necg \beta'\right) \wedge \necg \alpha < \necg\gamma\right).\]
  So it suffices to show that each disjunct of the consequent shows (ii).
  For the disjunct when \(\beta = \alpha\) the result is trivial, and if \(\beta \neq \alpha\) then we can use that one of the conjuncts is \(\necg \beta < \necg \alpha\) together with \(T \vdash \necg \beta < \necg \alpha \wedge \necg \alpha < \necg \gamma \to \necg \beta < \necg \gamma\) to derive the required formula.
\end{proof}

\begin{theorem}
  \label{advancing-function-GL}
  We have that \(T \vdash \rho_x \to \bigvee_{y \mid x \sqsubseteq^+ y}\lambda_y\).
\end{theorem}
\begin{proof}
  We proceed by the conversely wellfounded order given by \(\sqsubset^+\), so let \(x \in W\) and assume that for any \(y \sqsupset^+ x\) we have that \(T \vdash \rho_y \to \bigvee_{z \mid y \sqsubseteq^+ z} \lambda_z\).
  First, note that by definition of \(\lambda_x\) we have that \(T \vdash \rho_x \wedge \left(\bigwedge_{y \mid x \sqsubset y} \posg \lambda_y\right) \to \lambda_x\), so \(T \vdash \rho_x \wedge \left(\bigwedge_{y \mid x \sqsubset y} \posg \lambda_y\right) \to \bigvee_{y \mid x \sqsubseteq^+ y}\lambda_y\).
  Then it suffices to show that \(T \vdash \rho_x \wedge \left(\bigvee_{y \mid x \sqsubset y} \necg \neg\lambda_y\right) \to \bigvee_{y \mid x \sqsubseteq^+ y}\lambda_y\).
  We proceed as follows:
  \begin{align*}
    &T \vdash \left(\bigvee_{y \mid x \sqsubset y} \necg \neg \lambda_y\right) \to \left(\bigvee_{y \mid x \sqsubset y} \necg \neg \lambda_y \wedge\bigwedge_{z \mid x \sqsubset z \neq y} \necg \neg \lambda_y < \necg \neg \lambda_z\right)
    &&\text{by Lemma~\ref{finding-minimum-GL}}, \\
    &T \vdash \rho_x \wedge \left(\bigvee_{y \mid x \sqsubset y} \necg \neg \lambda_y\right) \to \rho_x \wedge \left(\bigvee_{y \mid x \sqsubset y} \necg \neg \lambda_y \wedge\bigwedge_{z \mid x \sqsubset z \neq y} \necg \neg \lambda_y < \necg \neg \lambda_z\right)
    &&\text{by propositional reasoning}, \\
    &T \vdash \rho_x \wedge \left(\bigvee_{y \mid x \sqsubset y} \necg \neg \lambda_y\right) \to \left(\bigvee_{y \mid x \sqsubset y} \rho_x \wedge \necg \neg \lambda_y \wedge\bigwedge_{z \mid x \sqsubset z \neq y} \necg \neg \lambda_y < \necg \neg \lambda_z\right)
    &&\text{by propositional reasoning}, \\
    &T \vdash \rho_x \wedge \left(\bigvee_{y \mid x \sqsubset y} \necg \neg \lambda_y\right) \to \left(\bigvee_{y \mid x \sqsubset y} \rho_y\right)
    &&\text{by definition of \(\rho_y\)}, \\
    &T \vdash \rho_x \wedge \left(\bigvee_{y \mid x \sqsubset y} \necg \neg \lambda_y\right) \to \left(\bigvee_{y \mid x \sqsubseteq^+ y} \lambda_y\right)
    &&\text{by I.H.},
  \end{align*}
  as desired.
\end{proof}

Since \(\rho_0 = \top\) we have the following thanks to the previous theorem.

\begin{corollary}[S1]
  We have \(T \vdash \bigvee_{x \in W}\lambda_x\).
\end{corollary}

Then (S4) will be obtained thanks to the following lemma.

\begin{lemma}
  \label{completeness-of-rho-GL}
  For any \(x \in W\) we have that \(T \vdash \rho_x \to \necg \rho_x\).
\end{lemma}
\begin{proof}
  It is easy to show that \(\rho_x\) is a \(\Sigma_1\)-sentence, the result follows from the \(\Sigma_1\)-completeness of \(\necg\).
\end{proof}

\begin{corollary}[S4]\label{property-S4-GL}
  We have \(T \vdash \lambda_x \to \necg \bigvee_{y \mid x \sqsubseteq^+ y}\lambda_y\).
\end{corollary}
\begin{proof}
  We have the following reasoning.
  \begin{align*}
    &T \vdash \rho_x \to \bigvee_{y \mid x \sqsubseteq^+ y} \lambda_y
    &&\text{by Theorem~\ref{advancing-function-GL}}, \\
    &T \vdash \necg\rho_x \to \necg\bigvee_{y \mid x \sqsubseteq^+ y} \lambda_y
    &&\text{by \((\NEC[\necg])\) and \((\Kax[\necg])\)}, \\
    &T \vdash \rho_x \to \necg\bigvee_{y \mid x \sqsubseteq^+ y} \lambda_y
    &&\text{by Lemma~\ref{completeness-of-rho-GL}}, \\
    &T \vdash \lambda_x \to \necg\bigvee_{y \mid x \sqsubseteq^+ y} \lambda_y
    &&\text{as \(T \vdash \lambda_x \to \rho_x\) by definition},
  \end{align*}
  as desired.
\end{proof}

\section{Arithmetical completeness for \(\ER\)}\label{sec:ER}

In this section, we prove the following arithmetical completeness theorem for the logic $\ER$\begin{theorem}[Arithmetical completeness for \(\ER\)]\label{arithmetical-completeness-ER}
  Let \(T_0,T_1\) be sound arithmetical theories and \(\prv[0], \prv[1]\) be predicates such that
  \begin{enumerate}
    \item for \(i \leq 1\), \(\prv[i]\) is a sound Hilbert-Bernays provability predicate of \(T_i\) in \(T_0\),
    \item \(T_0 \subseteq T_1\),
    \item \(T_1 \vdash \nec_0 \alpha \to \alpha\) for any sentence \(\alpha\), i.e., \(T_1\) proves the soundness of the \(\prv[0]\) provability predicate.
  \end{enumerate}
Then, for any formula \(\phi\) such that \(\ER \nvdash \phi\) there is an arithmetical interpretation \((\cdot)^*_{\prv[0], \prv[1]}\) with \(T_0 \nvdash (\phi)^*_{\prv[0],\prv[1]}\).
\end{theorem}
In proving this, the differences between our approach and the usual Solovay-style arithmetical completeness proofs via Kripke models will become clearer.
We might leave out some routine steps in the following section and instead refer to their respective steps in the previous section.

\subsection{Proofs and refutations}

\begin{figure}
  \[
    \AxiomC{}
    \RightLabel{\(\ax\)}
    \UnaryInfC{\(p, \Gamma \Rightarrow_i p, \Delta\)}
    \DisplayProof
    \qquad
    \AxiomC{}
    \RightLabel{\(\botL\)}
    \UnaryInfC{\(\bot, \Gamma \Rightarrow_i \Delta\)}
    \DisplayProof
    \qquad
    \AxiomC{\( \Gamma \Rightarrow_i \Delta\)}
    \RightLabel{\(\botR\)}
    \UnaryInfC{\( \Gamma \Rightarrow_i \bot,\Delta\)}
    \DisplayProof
  \]

  \[
    \AxiomC{\(\Gamma \Rightarrow_i \phi, \Delta\)}
    \AxiomC{\(\psi, \Gamma \Rightarrow_i \Delta\)}
    \RightLabel{\(\toL\)}
    \BinaryInfC{\(\phi \to \psi, \Gamma \Rightarrow_i \Delta\)}
    \DisplayProof
    \qquad
    \AxiomC{\(\phi, \Gamma \Rightarrow_i \psi, \Delta\)}
    \RightLabel{\(\toL\)}
    \UnaryInfC{\(\Gamma \Rightarrow_i \phi \to \psi, \Delta\)}
    \DisplayProof
  \]

  \[
    \AxiomC{\(\necd_i \Sigma_i, \nec_{\overline{i}} \Sigma_{\overline{i}}, \nec_i \phi \Rightarrow_i \phi\)}
    \RightLabel{\(\modal[i]{\ER}\)}
    \UnaryInfC{\(\nec_i \Sigma_i, \nec_{\overline{i}} \Sigma_{\overline{i}}, \Gamma \Rightarrow_j \nec_i \phi, \Delta\)}
    \DisplayProof
    \qquad
    \AxiomC{\(\necd_0 \phi, \Gamma \Rightarrow_1 \Delta\)}
    \RightLabel{\((\mathrm{ER})\)}
    \UnaryInfC{\(\nec_0 \phi, \Gamma \Rightarrow_1 \Delta\)}
    \DisplayProof
  \]
  \caption{\(\ER\) sequent rules}\label{fig:ER-sequent-rules}
  \separator
\end{figure}

We start by recalling the sequent calculus for \(\ER\) from~\cite{bimodal-ulip}.
A sequent in this calculus is going to be a triple \((\Gamma, \Delta, i)\) where \(\Gamma, \Delta\) are finite multisets of formulas and \(i \in \set{0,1}\).
We will denote the sequent \((\Gamma,\Delta,i)\) as \(\Gamma \Rightarrow_i \Delta\). 
The rules for this sequent calculus can be found at Figure~\ref{fig:ER-sequent-rules}.
As usual, the displayed formula at the conclusion of \((\botR)\), \((\toL)\), \((\toR)\), \((\modal[i]{\ER})\) is called the \emph{principal formula} of the rule instance.
The formulas of \(\nec_0 \Sigma_0, \nec_1 \Sigma\) at \((\modal[i]{\ER})\) are called the \emph{auxiliary formulas} of the rule instance.
Together, principal and auxiliary formulas are called the \emph{active formulas} of the instance.
Additionally, in \((\ax)\), \((\botL)\) and \((\modal[i]{\ER})\) the formulas at \(\Gamma\) and \(\Delta\) are called \emph{weakening formulas}.
The sequents which are instances of \((\ax)\) and \((\botL)\) are called \emph{axiomatic sequents}.
An \emph{extended axiomatic sequent} is a sequent \(\Gamma \Rightarrow_i \Delta\) such that \(\Gamma \cap \Delta \neq \varnothing\) or \(\bot \in \Gamma\).

\begin{definition}
  We define the sequent calculus \(\g{\ER}\) as the one given by the rules of Figure~\ref{fig:ER-sequent-rules}.
\end{definition}

It is straightforward to show that extended axiomatic sequents are provable.
\begin{proposition}
  Let \(S\) be an extended axiomatic sequent, then \(\g{\ER} \vdash S\).
\end{proposition}

To establish the correspondence of \(\g{\ER}\) with \(\ER\) we first define a translation from sequents to formulas as follows:
\[
  (\Gamma \Rightarrow_0 \Delta)^\flat = \bigwedge \Gamma \to \bigvee \Delta,
  \qquad\qquad
  (\Gamma \Rightarrow_1 \Delta)^\flat = \nec_1\left(\bigwedge \Gamma \to \bigvee \Delta\right).
\]

\begin{theorem}[\cite{bimodal-ulip}]
  \label{calculus-ER-sound-and-complete}
  For any sequent \(S\) we have that \(\ER \vdash S^\flat\) if and only if \(\g{\ER} \vdash S\).
\end{theorem}

\begin{figure}
  \[
    \AxiomC{\(\Gamma \Rightarrow_i \Delta\)}
    \RightLabel{\(\wk\)}
    \UnaryInfC{\(\Gamma, \Gamma' \Rightarrow_i \Delta, \Delta'\)}
    \DisplayProof
    \quad
    \AxiomC{\(\phi, \phi, \Gamma \Rightarrow_i \Delta\)}
    \RightLabel{\(\ctr\)}
    \UnaryInfC{\(\phi, \Gamma \Rightarrow_i \Delta\)}
    \DisplayProof
    \quad
    \AxiomC{\(\Gamma \Rightarrow_i \phi, \phi, \Delta\)}
    \RightLabel{\(\ctr\)}
    \UnaryInfC{\(\Gamma \Rightarrow_i \phi, \Delta\)}
    \DisplayProof
  \]
  \caption{\(\ER\) structural rules}
  \label{fig:ER-structural-rules}
  \separator
\end{figure}

\(\g{\ER}\) is also structurally well-behaved as the following lemma recalls.
\begin{lemma}[\cite{bimodal-ulip}]\label{lem:Wk-Ctr-Invert-ER}
  We have the following:
  \begin{enumerate}
    \item \((\wk)\) is height-preserving admissible in \(\g{\ER}\),
    \item \((\botR)\), \((\toL)\), \((\toR)\) and \((\mathrm{ER})\) are height-preserving invertible in \(\g{\ER}\),
    \item \((\ctr)\) is height-preserving admissible in \(\g{\ER}\).
  \end{enumerate}
\end{lemma}

\begin{figure}
  \[
    \AxiomC{\(\phi \to \psi, \Gamma \Rightarrow_i \phi, \Delta\)}
    \RightLabel{\(\toLsat_0\)}
    \UnaryInfC{\(\phi \to \psi, \Gamma \Rightarrow_i \Delta\)}
    \DisplayProof
    \qquad
    \AxiomC{\(\psi, \phi \to \psi, \Gamma \Rightarrow_i \Delta\)}
    \RightLabel{\(\toLsat_1\)}
    \UnaryInfC{\(\phi \to \psi, \Gamma \Rightarrow_i \Delta\)}
    \DisplayProof
    \qquad
    \AxiomC{\(\phi, \Gamma \Rightarrow_i \psi, \phi \to \psi,\Delta\)}
    \RightLabel{\(\toRsat\)}
    \UnaryInfC{\(\Gamma \Rightarrow_i \phi \to \psi, \Delta\)}
    \DisplayProof
  \]
  \caption{Proof search rules for \(\ER\)}
  \label{fig:ER-refutation-rules}
  \vspace{-0.3cm}
  \separator
\end{figure}

We define a notion of \emph{saturation}, similar to that in Definition~\ref{def:saturation-GL}, which will help us to exhaust the full proof search space. 

\begin{definition}
  Given a sequent \(S = (\Gamma \Rightarrow_i \Delta)\) we say that
  \begin{itemize}
    \item \(\phi\) is \emph{left-saturated} in \(S\) if
    \begin{enumerate}
      \item \(\phi\) is atomic or a \(\nec_1\)-formula,
      \item \(\phi = \phi_0 \to \phi_1\) and \(\phi_0 \in \Delta\) or \(\phi_1 \in \Gamma\),
      \item \(\phi = \nec_0 \phi_0\) and \(i = 0\) or \(\phi_0 \in \Gamma\).
    \end{enumerate} 
    \item \(\phi\) is \emph{right-saturated} in \(S\) if
    \begin{enumerate}
      \item \(\phi\) is atomic or a \(\nec\)-formula,
      \item \(\phi = \phi_0 \to \phi_1\), \(\phi_0 \in \Gamma\) and \(\phi_1 \in \Delta\).
    \end{enumerate} 
    \item \(S\) is said to be \emph{saturated} if every \(\phi \in \Gamma\) is left-saturated in \(S\) and every \(\psi \in \Delta\) is right-saturated in \(S\).
  \end{itemize}
  We define the saturation complexity of a sequent \(S = (\Gamma \Rightarrow_i \Delta)\) as
  \[
    \satc{S} = \sum \set{|\phi| \mid (\phi \in \Gamma \text{ and } \phi \text{ not left-saturated in }S) \text{ or }(\phi \in \Delta \text{ and } \phi \text{ not right-saturated in }S)},
  \]
  where \(|\phi|\) is the complexity of the formula.
  In words, the saturation complexity of a sequent is the sum of the complexities of the formulas in \(S\) which are not saturated.
\end{definition}

Finally, with the notion of saturation, we can define the desired game.

\begin{definition}[Proof Search Game]
  Given a sequent \(S\), we define the rule set for proof search in \(\ER\) denoted \(\ps_{\ER}\) as follows.
  \begin{enumerate}
    \item The positions are sequents \(S\).
    \item We will say that \(S'\) is a possible move at \(S\) if
      \begin{enumerate}
        \item \(S\) is not an extended axiomatic sequent, and
        \item \((S',S)\) is an instance of one of the rules \((\toLsat_0), (\toLsat_1), (\toRsat), (\mathrm{ER})\) with a non-saturated principal formula, or \(S\) is saturated and \((S',S)\) is a maximal instance of \((\modal[i]{\ER})\) for \(i \leq 1\).
      \end{enumerate}
    \item The positions of player \(\pI\) are the sequents which are saturated and not extended axiomatic.
      Notice that then, positions of player \(\pII\) are the sequents which are either non-saturated or extended axiomatic.
  \end{enumerate}
  It will be common to call player \(\pI\) \emph{verifier} and player \(\pII\) \emph{refuter}.
\end{definition}
We are going to show that for any sequent \(S\) the game \(\ps_{\ER} @ S\) is a short game.
We define a measure on sequents that is similar to the measure for \(\GL\):
\begin{align*}
  &|\Gamma \Rightarrow \Delta|_{\nec} = |\sub(\Gamma \Rightarrow \Delta)\setminus \set{\nec_i \phi \in \Gamma \mid i \leq 1}|, \\
  &|S|_{\ER} = \omega |S|_{\nec} + \satc{S}.
\end{align*}
\begin{lemma}
  If \(S,S'\) are sequents such that \(S'\) is a possible move at \(S\) then \(|S|_{\ER} > |S'|_{\ER}\).
\end{lemma}
\begin{proof}
  Analogous to the proof of Lemma~\ref{GL-property-of-the-measure}.
\end{proof}

As we did for \(\GL\), using the measure \(|\cdot|_{\ER}\) is then easy to establish that \(\ps_{\ER} @ S\) is a short game.

\begin{theorem}\label{short-game-psER}
  For any sequent \(S\), \(\ps_{\ER}@ S\) is a short game.
  As a corollary, either verifier or refuter have a winning strategy in \(\ps_{\ER} @ S\).
\end{theorem}
\begin{proof}
  Analogous to the proof of Theorem~\ref{short-game-psGL}
\end{proof}

As in the case of \(\GL\) we can show the correspondence between verifier's winning strategies and the existence of proofs.

\begin{theorem}
  \label{fundamental-theorem-of-refutations-ER}
  For any sequent \(S\) we have that
  \begin{enumerate}
    \item \(\g{\ER} \vdash S\) if and only if verifier has a winning strategy in \(\ps_{\ER} @ S\).
    \item \(\g{\ER} \not\vdash S\) if and only if refuter has a winning strategy in \(\ps_{\ER} @ S\), i.e., there is a refutation for \(S\).
  \end{enumerate}
\end{theorem}
\begin{proof}
  Analogous to Theorem~\ref{fundamental-theorem-of-refutations-GL}
\end{proof}

\subsection{Arithmetical completeness for \(\ER\)}\label{sec:arithmetical-completeness-ER}

For the remainder of Subsection~\ref{sec:arithmetical-completeness-ER} let us fix two arithmetical theories \(T_0,T_1\) and predicates \(\prv[0], \prv[1]\) satisfying the conditions displayed at Theorem~\ref{arithmetical-completeness-ER}.

Now, assume we have a formula \(\phi\) such that \(\ER \not\vdash \phi\), i.e., \(\g{\ER} \not \vdash { \Rightarrow_0 \phi}\).
This implies that there is a refutation of \(\Rightarrow_0 \phi\), let us denote this winning strategy for refuter by \(s_{\pII}\).
Since \(\ps_{\ER} @ ({\Rightarrow \phi})\) is a short game we have that \(\set{w \in \run({\Rightarrow_0 \phi}| s_{\pII}) \mid \last{w} \in P_{\pI}}\) of \(s_{\pII}\)-compatible runs whose last sequent is a position of verifier is also finite, let it be \(\set{w_1, \ldots, w_N}\).
Additionally, assume that we have picked the indexes such that if the run \(w_y\) extends the run \(w_x\) then \(x < y\) (this can be done easily because the set of runs with \(\leadsto^+\) is a finite tree).
Define \(W = \set{0,\ldots,N} \subseteq \mathbb{N}\), and define the relation \(x \sqsubset y\) for \(x,y \in W\) as
\begin{enumerate}
  \item \(x = 0 \) and \(y = 1\), or
  \item \(x,y > 0\), \(w_x = (S_0,\ldots,S_n)\), \(w_y = (S_0,\ldots,S_n,S_{n+1}, \ldots,S_m)\) for \(n < m\) and there is no position of verifier in \(S_{n+1}, \ldots,S_{m-1}\).
\end{enumerate}
For \(x \in W \setminus \set{0}\) let us denote the last sequent of \(w_x\) as \(S_x = \Gamma_x \Rightarrow_{i_x} \Delta_x\) and let us define \(S_0 := \Gamma_0 \Rightarrow_0 \Delta_0 := {\Rightarrow_0 \phi}\).
Since if \(w_y\) extends \(w_x\) we have that \(x < y\) we can see that the only position of verifier in \(w_1\) is the last one.
This means that the sequent \(S_1\) is obtained from \(S_0\) by just applying the rules \((\toL_0^{\mathrm{sat}})\), \((\toL_1^{\mathrm{sat}})\) and \((\toR^{\mathrm{sat}})\) so \(\Delta_0 \subseteq \Delta_1\).
This property will be useful later.
Finally we define the sets \(W_0 = \set{x \in W \setminus \set{0} \mid i_x = 0}\) and \(W_1 = \set{x \in W \setminus \set{0} \mid i_x = 1}\).

The same remarks that we made for the construction in $\GL$ at page~\pageref{remarks-about-GL-tree} could be made here.
We can view the structure $(W , \sqsubset)$ as a Kripke structure. However, just as transitivity is only enforced syntactically by how the sequent rules are defined, the ``pseudo reflexivity'' of the $\ER$ axiom gets solely enforced by the (ER) rule. Thus it is insufficient to consider only the underlying frame of these structures, which is how these refutations differ from e.g.\ Visser's models~\cite{Visser1995-VISACO-3}.
It also differs from the models used by Carlson~\cite{Carlson1986ModalLW} as they are reflexive and monotone (i.e., \(W_1 \subseteq W_0\) in Carlson's models while for us \(W_0 \cap W_1 = \varnothing\)).
In fact, the Kripke model is not even a model of \(\ER\), it only satisfies the necessary instances of the axioms of \(\ER\) to witness the unprovability of \(\phi\).
In the next subsubsection, we will see how this requires us to adapt some proofs to our setting.

The following lemma shows that, despite lacking the usual structural properties, the tree carries all the information necessary for the refutation of \(\phi\) in \(\ER\) by virtue of the shape of the rules of \(\g{\ER}\).
\begin{lemma}\label{properties-of-refutations-ER}
  We have the following for \(x \in W \setminus \set{0}\),
  \begin{enumerate}
    \item \(\bot \not \in \Gamma_x\), \(\Gamma_x \cap \Delta_x = \varnothing\) and \(\Gamma_x \Rightarrow \Delta_x\) is saturated.
    \item If \(\nec_i \psi \in \Delta_x\) then there is a \(y\) such that \(x \sqsubset y \in W_i\) and \(\psi \in \Delta_y\).
    \item If \(\nec_i \psi \in \Gamma_x\) and \(x \sqsubset^+ y\) then \(\nec_i \psi \in \Gamma_y\).
    Additionally, if \(y \in W_i\) then \(\psi \in \Gamma_y\).
    \item If \(\nec_0 \psi \in \Gamma_x\) and \(x \in W_1\) then \(\psi \in \Gamma_x\).
    \item If \(\nec_0 \psi \in \Gamma_x\) and \(x \sqsubset^+ y\) imply \(\psi \in \Gamma_y\).
  \end{enumerate}
\end{lemma}
\begin{proof}
  The proofs of the first three properties are analogous to their equivalents in Lemma~\ref{properties-of-refutations-GL}.
  The fourth property is a trivial consequence of the definition of saturation and the fifth property follows from properties three and four.
\end{proof}
We will use these properties of refutations in the remainder of the Subsection~\ref{sec:arithmetical-completeness-ER}, sometimes without explicitly mentioning them.
With these preliminaries in place, we are ready to prove the desired result.

\subsubsection{The heart of the proof}
In this subsection we show how to construct the arithmetical interpretation \((\cdot)^*\) such that \(T_0 \not \vdash \phi^*\).
For this we will assume the existence of some arithmetical sentences, called \emph{Solovay sentences}, \((\lambda_0,\ldots,\lambda_{N-1})\) such that
\begin{multicols}{2}
\begin{enumerate}\label{solovay-properties-ER}
  \item (\(\mathrm{S1}\)) \(T_0 \vdash \bigvee_{x \in W} \lambda_x\),
  \item (\(\mathrm{S2}\)) \(T_0 \vdash \neg(\lambda_x \wedge \lambda_y)\) for \(x \neq y\),
  \item (\(\mathrm{S3}\)) \(T_0 \vdash \lambda_x \to \posg_i \lambda_y\) for \(x \sqsubset y \in W_i\),
  \item (\(\mathrm{S4}\)) \(T_0 \vdash \lambda_x \to \necg_i\bigvee_{y \mid x \sqsubseteq^+ y} \lambda_y\),
  \item (\(\mathrm{S5}\)) \(T_0 \vdash \lambda_x \to \necg_i\neg \lambda_x\) for \(x \in W_i\)
\end{enumerate}
\end{multicols}
These conditions are the natural generalization of conditions (\(\mathrm{S1}\))-(\(\mathrm{S5}\)) for \(\GL\) at page~\pageref{solovay-properties-GL}.
The main difference is that now we have to index the modalities as well as respectively the set $W$ which corresponds to the index of the final sequent of $s_{\pII}$-compatible runs.
As we did for \(\GL\), in the next subsubsection, we will show how to construct such sentences using the ideas of Solovay.
We can also show that \(\lambda_0\) holds in the standard model.
\begin{lemma}
  We have that \(\mathbb{N} \vDash \lambda_0\).
\end{lemma}
\begin{proof}
  Since \(T_0 \vdash \bigvee_{x \in W} \lambda_x\) and \(T_0\) is sound we have that \(\mathbb{N} \vDash \bigvee_{x \in W} \lambda_x\).
  Then \(\mathbb{N} \vDash \lambda_x\) for some \(x \in W\), assume further that \(x > 0\) and derive a contradiction.
  Since \(x > 0\) we have that \(x \in W_i\) for some \(i \leq 1\) and then, using (S5) and soundness of \(T_0\), \(\mathbb{N}\vDash \lambda_x \to \necg_i \neg \lambda_x\), so \(\mathbb{N} \vDash \necg_i \neg \lambda_x\).
  As \(\prv[i]\) is a provability predicate for \(T_i\) we have that \(T_i \vdash \neg \lambda_x\) and since both theories are sound we conclude that \(\mathbb{N} \vDash \neg \lambda_x\), i.e., \(\mathbb{N} \not\vDash \lambda_x\), a contradiction.
\end{proof}

Additionally we have the following two lemmas.
\begin{lemma}\label{lm:er-prop}
  \(T_1 \vdash \neg \lambda_x\) for \(x \in W_0\).
\end{lemma}
\begin{proof}
  We have the following reasoning.
  \begin{align*}
    &T_0 \vdash \lambda_x \to \necg_0 \neg \lambda_x,
    &&\text{by \((\mathrm{S5})\)}; \\
    &T_1 \vdash \lambda_x \to \necg_0\neg \lambda_x,
    &&\text{since \(T_0 \subseteq T_1\)}; \\
    &T_1 \vdash \lambda_x \to \neg \lambda_x,
    &&\text{since \(T_1 \vdash \necg_0 \neg \lambda_x \to \neg \lambda_x\)}; \\
    &T_1 \vdash \neg \lambda_x ,
    &&\text{by logical reasoning}.
  \end{align*}
  As desired.
\end{proof}

\begin{lemma}\label{lm:mono}
  \(T_0 \vdash \necg_0 \alpha \to \necg_1 \alpha\) for any arithmetical sentence \(\alpha\).
\end{lemma}
\begin{proof}
  We have the following reasoning.
  \begin{align*}
    &T_1 \vdash \necg_0 \alpha \to \alpha,
    &&\text{assumption};\\
    &T_0 \vdash \necg_1\necg_0 \alpha \to \necg_1\alpha,
    &&\text{by \((\NEC[\necg])\) and \((\Kax[\necg])\) of \(\prv[1]\) in \(T_0\)};\\
    &T_0 \vdash \necg_0 \alpha \to \necg_1\alpha,
    &&\text{by \(\Sigma_1\)-completeness of \(\prv[1]\) in \(T_0\)}.
  \end{align*}
  As desired.
\end{proof}

We define the arithmetical interpretation
\[
  p^* := \bigvee_{x \mid p \in \Gamma_x} \lambda_x.
\]
This will be the desired interpretation such that \(T \nvdash \phi_0^*\).
As before, we need to establish a truth lemma to prove this.

\begin{lemma}[Truth Lemma]
  For any \(x \in W \setminus \set{0}\), we have that
  \begin{align*}
    &\psi \in \Gamma_x \text{ implies } T_0 \vdash \lambda_x \to \psi^*, 
    &&\psi \in \Delta_x \text{ implies } T_0 \vdash \lambda_x \to \neg \psi^*.
  \end{align*}
\end{lemma}
\begin{proof}
  We proceed by induction on \(\psi\).
  The cases in which \(\psi\) is a propositional variable, \(\bot\), or an implication are proved as in the case of \(\GL\) (see Lemma~\ref{truth-lemma-GL}).
  We are left with two cases, one for \(\nec_0\) and one for \(\nec_1\).

  Case \(\psi = \nec_0 \psi_0\), assume that \(\nec_0 \psi_0 \in \Gamma_x\).
  We further assume that \(x \in W_0\), then we have the following reasoning.
  \begin{align*}
    &T_0 \vdash \lambda_y \to \psi^*_0 \text{ for \(x \sqsubset^+ y\)}
    &&\text{by the I.H.\ since \(x \sqsubset^+ y\) and \(\nec_0 \psi_0 \in \Gamma_x\) implies \(\psi_0 \in \Gamma_y\);} \\
    &T_0 \vdash \left(\bigvee_{y \mid x \sqsubset^+ y} \lambda_y\right) \to \psi^*_0 
    &&\text{by logical reasoning;} \\
    &T_0 \vdash \necg_0\left(\bigvee_{y \mid x \sqsubset^+ y } \lambda_y\right) \to \necg_0 \psi_0^*
    &&\text{using \((\NEC[\necg])\) and \((\Kax[\necg])\) of \(\prv[0]\) in \(T_0\)}; \\
    &T_0 \vdash \necg_0\left(\bigvee_{y \mid x \sqsubseteq^+ y } \lambda_y\right) \to \necg_0 \psi_0^*
    &&\text{using (\(\mathrm{S}5\)) since \(x \in W_0\)}; \\
    &T_0 \vdash \lambda_x \to (\nec_0\psi_0)^* 
    &&\text{by (\(\mathrm{S}4\))}.
  \end{align*}
  Next, assume that \(x \in W_1\).
  Since \(\nec_0 \psi \in \Gamma_x\) we have that \(\psi \in \Gamma_x\).
  We have the following reasoning.
  \begin{align*}
    &T_0 \vdash \lambda_y \to \psi^*_0 \text{ for \(x \sqsubset^+ y\)}
    &&\text{by the I.H.\ since \(x \sqsubset^+ y\) and \(\nec_0 \psi_0 \in \Gamma_x\) implies \(\psi_0 \in \Gamma_y\);} \\
    &T_0 \vdash \lambda_x \to \psi^*_0 
    &&\text{by the I.H.\ since \(x \in W_1\) and \(\nec_0 \psi_0 \in \Gamma_x\) implies \(\psi_0 \in \Gamma_x\);} \\
    &T_0 \vdash \left(\bigvee_{y \mid x \sqsubseteq^+ y} \lambda_y\right) \to \psi^*_0 
    &&\text{by logical reasoning from the two previous lines;} \\
    &T_0 \vdash \necg_0\left(\bigvee_{y \mid x \sqsubseteq^+ y } \lambda_y\right) \to \necg_0\psi_0^* 
    &&\text{using \((\NEC[\necg])\) and \((\Kax[\necg])\) of \(\prv[0]\) in \(T_0\)}; \\
    &T_0 \vdash \lambda_x \to (\nec_0\psi_0)^* 
    &&\text{by (\(\mathrm{S}4\))}.
  \end{align*}
  Finally, suppose that \(\nec_0 \psi_0 \in \Delta_x\).
  We know that there is a \(y\) such that \(x \sqsubset y \in W_0\) and \(\psi_0 \in \Delta_y\).
  We have the following reasoning.
  \begin{align*}
    &T_0 \vdash \lambda_y \to \neg \psi_0^*,
    &&\text{by I.H.;} \\
    &T_0 \vdash \psi_0^* \to \neg \lambda_y,
    &&\text{by logical reasoning;} \\
    &T_0 \vdash \necg_0 \psi_0^* \to \necg_0\neg \lambda_y,
    &&\text{using \((\NEC[\necg])\) and \((\Kax[\necg])\) of \(\prv[0]\) in \(T_0\)}; \\
    &T_0 \vdash \posg_0 \lambda_y \to \neg(\nec_0 \psi_0)^*,
    &&\text{by logical reasoning;} \\
    &T_0 \vdash \lambda_x \to \neg (\nec_0 \psi_0)^*,
    &&\text{using \((\mathrm{S3})\) since \(x \sqsubset y \in W_0\).}
  \end{align*}
  This establishes the desired result.

  Case \(\psi = \nec_1 \psi_0\), assume that \(\nec_1 \psi_0 \in \Gamma_x\).
  First, assume that \(x \in W_0\), then we have the following reasoning.
  \begin{align*}
    &T_0 \vdash \lambda_y \to \psi^*_0 \text{ for \(x \sqsubset^+ y \in W_1\)}
    &&\text{by the I.H.\ since \(x \sqsubset^+ y \in W_1\) and \(\nec_1 \psi_0 \in \Gamma_x\) implies \(\psi_0 \in \Gamma_y\);} \\
    &T_0 \vdash \left(\bigvee_{y \mid x \sqsubset^+ y \in W_1} \lambda_y\right)\to \psi^*_0 
    &&\text{by logical reasoning;} \\
    &T_0 \vdash \necg_1\left(\bigvee_{y \mid x \sqsubset^+ y \in W_1} \lambda_y\right) \to \necg_1\psi_0^* 
    &&\text{using \(T_0 \subseteq T_1\) and \((\NEC[\necg])\), \((\Kax[\necg])\) of \(\prv[1]\) in \(T_0\)}; \\
    &T_0 \vdash \necg_1\left(\bigvee_{y \mid x \sqsubset^+ y} \lambda_y\right) \to (\nec_1\psi_0)^*
    &&\text{using that \(T_0 \vdash \necg_1\left(\bigwedge_{y \in W_0} \neg \lambda_y\right)\) by Lemma~\ref{lm:er-prop}},\\
    &T_0 \vdash \necg_1\left(\bigvee_{y \mid x \sqsubseteq^+ y} \lambda_y\right) \to (\nec_1\psi_0)^*
    &&\text{by  (\(\mathrm{S}5\)) since \(x \in W_0\) and \(T_0 \vdash \necg_0\neg \lambda_x \to \necg_1\neg \lambda_x\)};\\
    &T_0 \vdash \lambda_x \to (\nec_1\psi_0)^* 
    &&\text{by \((\mathrm{S4})\) and the two previous lines}.
  \end{align*}
  Now, assume that \(x \in W_1\), we have the following reasoning.
  \begin{align*}
    &T_0 \vdash \lambda_y \to \psi^*_0 \text{ for \(x \sqsubset^+ y \in W_1\)},
    &&\text{by the I.H.\ since \(x \sqsubset^+ y \in W_1\) and \(\nec_1 \psi_0 \in \Gamma_x\) implies \(\psi_0 \in \Gamma_y\);} \\
    &T_0 \vdash \left(\bigvee_{y \mid x \sqsubset^+ y \in W_1} \lambda_y\right) \to \psi^*_0,
    &&\text{by logical reasoning;} \\
    &T_0 \vdash \necg_1\left(\bigvee_{y \mid x \sqsubset^+ y \in W_1} \lambda_y\right) \to \necg_1\psi_0^*,
    &&\text{by \(T_0 \subseteq T_1\) and \((\NEC[\necg])\), \((\Kax[\necg])\) of \(\prv[1]\) in \(T_0\)}; \\
    &T_0 \vdash \necg_1\left(\bigvee_{y \mid x \sqsubset^+ y} \lambda_y\right) \to (\nec_1\psi_0)^*,
    &&\text{using that \(T_0 \vdash \necg_1\left(\bigwedge_{y \in W_0} \neg \lambda_y\right)\) by Lemma~\ref{lm:er-prop}};\\
    &T_0 \vdash \necg_1\left(\bigvee_{y \mid x \sqsubseteq^+ y} \lambda_y\right) \to (\nec_1\psi_0)^*,
    &&\text{by \((\mathrm{S}5)\) since \(x \in W_1\)};\\
    &T_0 \vdash \lambda_x \to (\nec_1\psi_0)^*,
    &&\text{by (\(\mathrm{S}4\))}.
  \end{align*}
  Finally, assume that \(\nec_1 \psi_0 \in \Delta_x\).
  We know that there is a \(y\) such that \(x \sqsubset y \in W_1\) and \(\psi_0 \in \Delta_y\).
  We have the following reasoning.
  \begin{align*}
    &T_0 \vdash \lambda_y \to \neg \psi_0^*,
    &&\text{by I.H.;} \\
    &T_0 \vdash \psi_0^* \to \neg \lambda_y,
    &&\text{by logical reasoning;} \\
    &T_0 \vdash \necg_1\psi_0^* \to \necg_1\neg \lambda_y,
    &&\text{using \(T_0 \subseteq T_1\) and \((\NEC[\necg])\), \((\Kax[\necg])\) of \(\prv[1]\) in \(T_0\)}; \\
    &T_0 \vdash \posg_1 \lambda_y \to \neg(\nec_1\psi_0)^*
    &&\text{by logical reasoning;} \\
    &T_0 \vdash \lambda_x \to \neg (\nec_1 \psi_0)^*
    &&\text{using \((\mathrm{S3})\) since \(x \sqsubset y \in W_1\).}
  \end{align*}
  We conclude the desired result.
\end{proof}




The proof of completeness proceeds as follows.
Assume that \(T_0 \vdash \phi^*\) and let us reach a contradiction.
We know that \(\set{\phi} = \Delta_0 \subseteq \Delta_1\) and \(0 \sqsubset 1\), we also assume that \(1 \in W_i\).
Using these facts we have the following proof:
\begin{align*}
  &T_0 \vdash \lambda_1 \to \neg \phi^*,
  &&\text{by Truth Lemma since \(\phi \in \Delta_1\)}; \\
  &T_0 \vdash \neg \lambda_1,
  &&\text{since \(T_0 \vdash \phi^*\)}; \\
  &T_0 \vdash \necg_i\neg \lambda_1,
  &&\text{by \((\NEC[\necg])\)}; \\
  &T_0 \vdash \neg \lambda_0
  &&\text{by \((\mathrm{S3})\) as \(0 \sqsubset 1 \in W_i\)}; \\
  &\mathbb{N} \not\vDash \lambda_0,
  &&\text{by soundness of \(T_0\)}.
\end{align*}
This yields the desired contradiction, since we already established that \(\mathbb{N} \vDash \lambda_0\).

\begin{example}
  Let us remember from Example~\ref{example-GL} that the theories \(\EA\), \(\PA\) and \(\IS{n}\) for \(n \geq 1\) have natural provability predicates \(\prv[\EA]\), \(\prv[\PA]\) and \(\prv[\IS{n}]\) respectively.
  Then, if we consider these theories with their natural provability predicates, the classical examples of arithmetical completeness for \(\ER\) (\cite{Carlson1986ModalLW, Beklemishev1996-BEKBLF}), owing to the essential reflexivity of \(\PA\),  include \((\EA, \PA)\), \((\IS{n},\PA)\) for \(n < \omega\). If we extend the result to include set theories we would have that \(\ER\) is also arithmetically complete with respect to \((\PA, \mathsf{ZF})\).
\end{example}

\subsubsection{Construction of the Solovay function}

We are going to define the sentences \(\lambda_0, \ldots, \lambda_N\) fulfilling the properties (S1) to (S5) at page~\pageref{solovay-properties-ER}.
The idea is the same as for \(\GL\), the only difference being that in order to move to a node in \(W_i\) we need a proof in the corresponding \(\prv[i]\).
Again \(\lambda_x\) will claim that the limit of the Solovay function is \(x\) and \(\rho_x\) will claim that the Solovay function passes through \(x\).
First, we recursively define the arithmetical sentences \(\rho_x(p_0,\ldots,p_N)\) as
\begin{align*}
  &\rho_0(p_0,\ldots,p_N) := \top, && \\
  &\rho_x(p_0,\ldots,p_N) := \rho_y(p_0,\ldots,p_N) \wedge \prv[i_x](\dot\neg p_x) \wedge \bigwedge_{z \mid y \sqsubset z \neq x} \prv[i_x](\dot\neg p_x) < \prv[i_z]( \dot\neg p_z)
  &&\text{where }y \sqsubset x.
\end{align*}
In the definitions above, we define \(\prv[i](\dot{\neg} x)\) as a strict \(\Sigma_1\)-formula equivalent in \(T\) to the formula \(\exists p, y (\mathrm{Neg}(x,y) \wedge \prf[i](p,y))\), i.e., as saying that the negation of \(x\) is provable in \(\prv[i]\).
Note that such a strict \(\Sigma_1\)-formula exists since \(\mathrm{Neg}(x,y)\) is \(\Delta_0\) and then the formula above is \(\Sigma_1\).

The \(\rho_x\)'s are well-defined because each \(x \in W \setminus \set{0}\) has a unique \(y \in W\) such that \(y \sqsubset x\), since \(\sqsubset\) is the immediate successor relation of a tree.
We consider the system of formulas
\[
  \rho_x(p_0,\ldots,p_N) \wedge \bigwedge_{y \mid x \sqsubset y} \con[i_y](p_y)\qquad \text{for }x \in \set{0,\ldots,N}.
\]
where we define \(\con[i](x) := \neg \prv[i](\dot{\neg} x)\), i.e.,\ the consistency predicate of \(\prv[i]\).
Given an arithmetical sentence \(\alpha\) we will write \(\posg_i \alpha\) to mean \(\con[i](\Gdl{\alpha})\).

By the uniform diagonal lemma for \(T\), we obtain arithmetical sentences  \(\lambda_0,\ldots,\lambda_N\) such that
\[
  T_0 \vdash \lambda_x \leftrightarrow \rho_x \wedge \bigwedge_{y \mid x \sqsubset y} \posg_{i_y} \lambda_y, \qquad \text{for }x \in \set{0,\ldots,N}.
\]
where \(\rho_x\) is just \(\rho_x := \rho_x(\ulcorner \lambda_0\urcorner, \ldots,\ulcorner\lambda_N\urcorner)\) for \(x \in W\).
Additionally, by definition of the \(\rho_x\)'s we have
\[
  T_0 \vdash \rho_x \leftrightarrow \rho_y \wedge \necg_{i_x} \neg \lambda_x \wedge \bigwedge_{z \mid y \sqsubset z \neq x} \necg_{i_x} \neg \lambda_x < \necg_{i_z} \neg \lambda_z \qquad \text{where }y \sqsubset x;
\]
We also assume that, for \(x \neq y\), we have \(\lambda_x \neq \lambda_y\), i.e.,\ that for different elements of \(W\) the Solovay sentences are not syntactically equal.
Otherwise, if \(\lambda_x = \lambda_y\) for \(x \neq y\), we could redefine \(\lambda_y\) as \(\lambda_x \wedge \top\) and as there are only finitely many Solovay sentences this modification suffices.
As immediate consequences of these definitions we have the following propositions.

\begin{proposition}[S3]
  We have that \(T_0 \vdash \lambda_x \to \posg_i \lambda_y\) for \(x \sqsubset y \in W_i\).
\end{proposition}

\begin{proposition}[S5]
  We have that \(T_0 \vdash \lambda_x \to \necg_i \neg \lambda_x\) for \(x \in W_i\).
\end{proposition}

As with \(\GL\) now we proceed to prove (S2).
We start by showing that if \(x\) and \(y\) do not belong to the same branch then the Solovay function cannot pass through both of them.

\begin{lemma}\label{incomparable-absurd-ER}
  If \(x \not\sqsubseteq^+ y\) and \(y \not \sqsubseteq^+ x\) then \(T_0 \vdash \neg(\rho_x \wedge \rho_y)\).
\end{lemma}
\begin{proof}
  Analogous to the proof of Lemma~\ref{incomparable-absurd-GL}.
\end{proof}

\begin{lemma}\label{limit-absurd-ER}
  If \(x \sqsubset y\) then \(T_0 \vdash \neg(\lambda_x \wedge \rho_y)\).
\end{lemma}
\begin{proof}
  Analogous to the proof of Lemma~\ref{limit-absurd-GL}.
\end{proof}

\begin{theorem}[S2]
  We have that \(T_0 \vdash \neg(\lambda_x \wedge \lambda_y)\) for \(x \neq y\).
\end{theorem}
\begin{proof}
  Analogous to the proof of Theorem~\ref{property-S2-GL}.
\end{proof}

The final property we need to establish, is that, if the Solovay function passes through a node \(x\), then the limit must be either \(x\) or a node above \(x\), from this all the rest of Solovay conditions can be derived.
To establish this, we first need the following lemma.

\begin{lemma}\label{finding-minimum-ER}
  We have that \(T_0 \vdash \left(\bigvee_{y \mid x \sqsubset y} \necg_{i_y} \neg \lambda_y\right) \to \left(\bigvee_{y \mid x \sqsubset y} \necg_{i_y} \neg \lambda_y \wedge \bigwedge_{z \mid x \sqsubset z \neq y} \necg_{i_y} \neg \lambda_y < \necg_{i_z} \neg \lambda_z\right)\).
\end{lemma}
\begin{proof}
  Analogous to the proof of Lemma~\ref{finding-minimum-GL}.
\end{proof}

\begin{theorem}\label{advancing-function-ER}
  We have that \(T_0 \vdash \rho_x \to \bigvee_{y \mid x \sqsubseteq^+ y}\lambda_y\).
\end{theorem}
\begin{proof}
  The proof is analogous to the proof of Theorem~\ref{advancing-function-GL}.
\end{proof}

Since \(\rho_0 = \top\), the preceding lemma yields the following corollary.

\begin{corollary}[S1]
  We have \(T_0 \vdash \bigvee_{x \in W}\lambda_x\).
\end{corollary}

\begin{lemma}\label{lem:completeness-of-rho-ER}
  For any \(x \in W\) and \(i \leq 1\) we have that \(T_0 \vdash \rho_x \to \necg_i \rho_x\).
\end{lemma}
\begin{proof}
  Analogous to the proof of Lemma~\ref{completeness-of-rho-GL}.
\end{proof}

\begin{corollary}[S4]
  We have \(T_0 \vdash \lambda_x \to \necg_i \bigvee_{y \mid x \sqsubseteq^+ y}\lambda_y\) for \(i \in \set{0,1}\).
\end{corollary}
\begin{proof}
  Analogous to the proof of Corollay~\ref{property-S4-GL}.
\end{proof}

\section{Future work}
\label{sec:future-work}

As future work, we would like to continue developing the syntactic approach for provability logic, investigating the uniformity of the method by applying it to other provability logics and also trying to expand the approach to provide new results.
Some natural directions for future research are as follows.
\begin{enumerate}
  \item Give a proof of the arithmetical completeness of \(\GLP\) using an analytic sequent calculus.
    Given the intricacies of \(\GLP\) we expect that it is not possible to obtain a cut-free pure sequent calculus, for a cut-free nested sequent calculus, see~\cite{shamkanov2015nested}).
  \item Give a proof of arithmetical completeness for interpretability logics, thus providing a more unified methodology that can treat both provability and interpretability (contrary to the semantic approach where some variant of Kripke semantics is needed).
  \item Give a proof of arithmetical completeness for intuitionistic modal logic, investigating if it is possible to simplify current proofs (which depend on characterizing the \(\Sigma_1\)-provability logic of \(\mathsf{HA}\) and then lift the result to the usual provability logic) with syntactic tools~\cite{mojtahedi2018TheSigma1, Mojtahedi2024, mojtahedi2026provabilitylogicha}.
  \item Kogure and Kurahashi are developing a program of finding arithmetical completeness of non-normal modal logic with respect to arithmetical predicates that do not satisfy all the usual conditions of provability, e.g., see~\cite{KOGURE2023103271, kogure2026modallogicsconjunctivelyclosed}.
  For this project they use different semantic frameworks such as a relational adaptation of neighborhood semantics (instead of Kripke models) or a new semantics based on closure predicates.
  Can we provide alternative proofs using the sequent calculus for these logics?
  This endeavour would also help us assess the flexibility and uniformity of our method.
\end{enumerate}

\section*{Acknowledgements}
The second author would like to thank Lev Beklemishev for his interest in this work and for drawing a parallel between Japaridze's ideas and the ideas presented in this article. Without his familiarity with Japaridze's work it would not have been possible for the authors to recognize the connection.

The second author is supported by the Swiss National Science Foundation project 200021\_214820.

\bibliography{bib}
\bibliographystyle{plain}
\end{document}